\documentclass[12pt]{article}
\usepackage{amssymb}
\usepackage{amsmath}
\usepackage{latexsym}
\newtheorem{theorem}{Theorem}[section]
\newtheorem{proposition}[theorem]{Proposition}
\newtheorem{corollary}[theorem]{Corollary}

\newtheorem{remark}[theorem]{Remark}

\begin{document}


\title{Genus Two Partition and Correlation Functions for Fermionic Vertex Operator Superalgebras~I}
\author{
 Michael P. Tuite and Alexander Zuevsky\thanks{
Supported by a Science Foundation Ireland Frontiers of Research Grant, and
by Max--Planck Institut f\"{u}r Mathematik, Bonn} \\
School of Mathematics,  
Statistics and Applied Mathematics,  
\\
National University of Ireland Galway \\
University Road, Galway, Ireland.\\
}
\maketitle
\begin{abstract}
We define the partition and $n$-point correlation functions for a vertex operator
superalgebra on a genus two Riemann surface formed by sewing two tori together.
For the free fermion vertex operator superalgebra we obtain a closed formula for the genus two continuous orbifold partition function in
 terms of an infinite dimensional determinant with entries arising from torus Szeg\"o kernels. 
We prove that the partition function is holomorphic in the sewing
parameters on a given suitable domain and describe its modular properties. Using the bosonized formalism, a new genus two Jacobi product identity is described for the Riemann theta series.
We compute and discuss the modular properties of the generating function for all $n$-point functions in terms of a genus two Szeg\"o kernel determinant. We also show that the
Virasoro vector one point function satisfies a genus two Ward identity.
\end{abstract}
\maketitle
\newpage
\section{Introduction}
Genus two (and higher) partition functions and
correlation functions have been studied for some time in string and conformal field theory e.g. 
\cite{EO}, \cite{FS}, \cite{DP}, \cite{Kn}, \cite{DVPFHLS}. Meanwhile, in the theory of Vertex Operator 
Algebras (VOAs) \cite{B}, \cite{FHL}, \cite{FLM}, \cite{Ka}, \cite{MN}, \cite{MT5} higher genus approaches based 
on algebraic geometry have also been developed e.g. \cite{TUY}, \cite{KNTY}, \cite{Z2}, \cite{U}. 
A more constructive VOA approach has recently been described whereby genus two partition and $n$-point 
correlation functions are defined in terms of genus one VOA data \cite{T}, \cite{MT1}, \cite{MT2}, \cite{MT3}, \cite{MT4}. 
This approach is based solely on the properties of a VOA with
no assumed analytic or modular properties for partition or correlation functions.
A compact genus two Riemann surface can be obtained from tori by 
either sewing two separate tori together, which we refer to as the 
$\epsilon$-formalism, or by self-sewing a torus, which we refer to as 
the $\rho$-formalism \cite{MT2}.  The theory of partition and 
$n$-point correlation functions in the $\epsilon$-formalism is 
described in ref. \cite{MT1} where these functions are explicitly 
computed for the Heisenberg VOA and its modules including lattice 
VOAs. The corresponding functions are considered in the 
$\rho$-formalism in ref. \cite{MT3}.  
 
 \medskip 
This paper extends these methods to the study of genus two partition and 
$n$-point functions in the $\epsilon$-formalism for Vertex Operator 
Superalgebras (VOSA). 
In particular, we explicitly compute and prove convergence 
and modular properties of the genus two continuous orbifold 
partition and $n$-point functions for the rank two fermion VOSA 
$V(H,\mathbb{Z}+\frac{1}{2})^{\otimes 2}$. (The alternative $\rho$-formalism is 
considered elsewhere \cite{TZ3}). These functions are computed in terms of 
appropriate torus $n$-point functions described in \cite{MTZ}. We also make 
extensive use of the expression of the genus two Szeg\"o  kernel $S^{(2)}$ of \eqref{eq:Szego2} in 
terms of genus one Szeg\"o kernel data described in \cite{TZ1}.   The 
partition function is then expressed as a certain infinite determinant 
whose components arise from genus one Szeg\"o kernel data. 
Furthermore, the generating function of all $n$-point correlation 
functions is computed in terms of a genus two Szeg\"o kernel 
determinant.

\medskip 
Section~2 consists of a review of aspects of the $\epsilon$-formalism for constructing a genus two Riemann surface by  sewing two separate tori with modular parameters $\tau_{1},\tau_{2}$ respectively 
for $(\tau_{1},\tau_{2},\epsilon)\in\cal{D}^\epsilon$, a specific domain for which the sewing 
is defined \cite{MT2}. We also review the construction of the genus two Szeg\"o kernel $S^{(2)}$ 
in terms of genus one Szeg\"o kernel data \cite{TZ1}. In particular we introduce an infinite block matrix
\begin{equation*}
Q=\left(
\begin{array}{cc}
0 & \xi F_1(\tau_{1})\\
-\xi F_2(\tau_{2}) & 0
\end{array}
\right),
\end{equation*}  
where $\xi=\pm\sqrt{-1}$ and $F_{a}(\tau_{a})$ for $a=1,2$ are certain infinite matrices whose entries 
involve twisted modular forms in $\tau_{a}$  associated with genus one Szeg\"o kernels \cite{MTZ}. 
Section~3 is a review of Vertex Operator Superalgebras (VOSA) and the Li--Zamolodchikov (Li--Z) metric 
on a VOSA \cite{L}, \cite{Sche}. The free fermion rank one VOSA $V(H,\mathbb{Z}+\frac{1}{2})$ is also reviewed. 
In Section~4 we consider the orbifold partition and $n$-point function on a genus two surface in the 
$\epsilon$-formalism for a VOSA with a Li--Z metric. 
These are defined in terms of genus one $n$-point orbifold functions associated with a pair of commuting 
VOSA automorphisms 
$f_a,g_a$ on a torus with modular parameter $\tau_a$ for $a=1,2$.
\medskip

Section~5 contains the main results of the paper wherein the partition function and the 
generating function for  $n$-point functions are computed for the rank two fermion VOSA with continuous 
automorphisms generated by the Heisenberg vector. In particular 
we prove in Theorem \ref{Theorem_Z2_boson} that the partition function is given by  
\begin{equation*}
Z^{(2)}\left[{f \atop g} \right](\tau _{1},\tau _{2}, \epsilon )=  
Z^{(1)}\left[{f_1 \atop g_1} \right](\tau_{1})\; 
Z^{(1)}\left[{f_2 \atop g_2} \right](\tau _{2}) \; \det
\left(I-Q\right),
\end{equation*}
where $f=(f_1,g_1)$ and $g=(f_2,g_2)$ and 
$Z^{(1)}\left[{f_a \atop g_a} \right](\tau_{a})$ 
is the orbifold partition function on the torus with modular parameter $\tau_a$. The partition function is 
holomorphic for $(\tau_{1},\tau_{2},\epsilon)\in\cal{D}^\epsilon$, a specific domain on which the $\epsilon$-formalism can be carried out \cite{MT2}. 
In Theorem~\ref{generating_n_point_rank_two_1} we find the generating 
function for all genus two $n$-point functions as a differential form
which is expressed in terms of a finite dimensional determinant of genus two Szeg\"o kernels $S^{(2)}$. 
We also discuss the bosonization of the fermion VOSA wherein the partition function can be expressed in 
terms of a genus two Riemann theta series and the Heisenberg genus two partition function. 
This leads to a new genus two version of the classical Jacobi product identity expressing the genus two 
Riemann theta series in terms of certain infinite products. 
We also discuss the genus two Ward identity satisfied by the Virasoro one point function in this bosonized setting.

\medskip
In Section~6 we discuss modular invariance of the genus two partition and $n$-point generating form under 
a modular group preserving $\cal{D}^\epsilon$.
The Appendix describes some general aspects of Riemann surfaces such as the period matrix, the projective 
connection and the prime form.
We also recall some facts from the classical and twisted  
elliptic function theory \cite{MTZ}.  

\medskip
We collect here notation for some of the more frequently occurring functions
and symbols employed. $\mathbb{Z}$ is the set of
integers, $\mathbb{C}$ the complex numbers, $\mathbb{H}$ the complex
upper-half plane. We will always take $\tau $ to lie in $\mathbb{H}$, and $z$
will lie in $\mathbb{C}$ unless otherwise noted. For a symbol $z$ we set $q_{z}=\exp (z)$ and in particular $q=q_{2 \pi i \tau}=\exp (2 \pi i \tau)$.  
\section{The Szeg\"o Kernel on a Genus Two Riemann Surface Formed from Two Sewn Tori}
\label{Sect_Epsilon_g}
The central role played by the Szeg\"o kernel $S^{(g)}$ for the fermion VOSA has been long known 
e.g. \cite{RS}, \cite{R}, \cite{DVFHLS}, \cite{DVPFHLS}. 
In this Section we review the form of the Szeg\"o kernel on a
Riemann surface $\Sigma^{(2)}$ of genus two obtained by sewing together two tori described in \cite{TZ1}. 
Some further details appear in Appendix~\ref{PrimeForm}.
\subsection{The Szeg\"o Kernel on a Riemann Surface}  
\label{Szegokernel}
Consider a compact connected Riemann surface $\Sigma^{(g)}$ of genus $g$ with canonical
homology cycle basis $a_{i},b_{i}$ for $i=1,\ldots, g$. Let  $\nu^{(g)}_{i}$ be a basis of holomorphic 1-forms 
with normalization $\oint_{a_{i}}\nu^{(g)}_{j}=2\pi i\delta_{ij}$ and
period matrix $\Omega^{(g)}_{ij}=\frac{1}{2\pi i}\oint_{b_{i}}\nu^{(g)}_{j}\in \mathbb{H}_g$, the Siegel upper 
half plane (e.g. \cite{FK}, \cite{Sp}).
Define the theta function with real characteristics  \cite{M}, \cite{F1},  \cite{FK}
\begin{eqnarray}
\label{theta}
&\vartheta^{(g)} \left[ {{\alpha} \atop {\beta} }\right] \left( {z} \vert    \Omega^{(g)}  \right)  
= 
\sum\limits_{ n \in {\mathbb{Z}^{g}} } 
e^{ i \pi (n + {\alpha}).\Omega^{(g)} .(n+{\alpha}) + 
 (n+ {\alpha}). ({{z}+ 2 \pi i {\beta}}) },&
\end{eqnarray}
for $\alpha=(\alpha_j),{\beta}=(\beta_j)\in \mathbb{R}^g$ and ${z}=(z_j)\in \mathbb{C}^g$ 
for $j=1,\ldots, g$.

The Szeg\"o Kernel \cite{Schi}, \cite{HS}, \cite{F1}, \cite{F2} 
is defined for $\vartheta \left[ {{\alpha} \atop {\beta}}\right] ( 0\vert \Omega^{(g)})\neq 0$ by 
\begin{equation}
\label{gSzego10}
 S^{(g)}\left[ {{\theta} \atop {\phi} }\right] (x, y) =
\frac
{ \vartheta^{(g)} \left[ {{\alpha} \atop {\beta}} \right] 
\left( \int_{y}^{x}\nu^{(g)}  \,\vert \Omega^{(g)}\right)} 
 {\vartheta^{(g)}\left[ {{\alpha} \atop {\beta}}\right] ( 0\vert \Omega^{(g)})  E^{(g)}(x, y)},
\end{equation}
where 
$\theta=({\theta}_j),\ \phi=(\phi_j)\in U(1)^n$ for
\begin{equation}
{\theta}_j=-e^{-2 \pi i \beta_j}, \quad{\phi}_j= -e^{2 \pi i \alpha_j},\quad j=1,\ldots, g, 
\label{eq:periodicities}
\end{equation}
and $E^{(g)}(x,y)$ is the prime form (see Appendix~\ref{PrimeForm}).
The factors of $-1$ in (\ref{eq:periodicities}) 
are included for later convenience.
The Szeg\"o kernel has multipliers along the $a_i$ and $b_j$ cycles in
 $x$ given by $-\phi_i$ and $-\theta_j$ respectively and is 
a meromorphic $(\frac{1}{2},\frac{1}{2})$-form satisfying  
\begin{eqnarray*}
  S^{(g)} \left[ {\theta \atop \phi} \right] (x,y) 
& \sim &  \frac{1}{x-y}\; dx^{\frac{1}{2}}\; dy^{\frac{1}{2}}
\quad \mbox{for } x\sim y, \label{Sz_local}\\
 S^{(g)}\left[ {\theta \atop \phi} \right] (x,y) 
& = &   -S^{(g)}\left[ {\theta^{-1} \atop \phi^{-1} } \right] (y,x),\label{Sz_skewsym}
\end{eqnarray*}
where $\theta^{-1}=({\theta}_i^{-1})$ and $\phi^{-1}=(\phi_i^{-1})$.  
\subsection{Genus Two Riemann Surfaces Formed from Two Sewn Tori}
\label{subsect_genus_two_eps}
Consider the genus two Riemann surface formed by sewing together two tori in the sewing scheme referred to 
as the $\epsilon$-formalism in refs.
\cite{MT1}, \cite{MT2}, \cite{TZ1}. Let $\Sigma_a^{(1)}=\mathbb{C}/{\Lambda}_{a}$ for $a=1,2$
be oriented tori with lattice 
${\Lambda}_{a}=2\pi i(\mathbb{Z}\tau _{a}\oplus \mathbb{Z})$ for 
$\tau _{a}\in \mathbb{H}$.  
Choose a local coordinate $z_{a}\in \mathbb{C}/{\Lambda}_{a}$ on $\Sigma^{(1)}_a$ in the
neighborhood of a point $p_{a}\in\Sigma^{(1)}_a$ and consider the closed disk $\left\vert
z_{a}\right\vert \leq r_{a}$ for $r_{a}<\frac{1}{2}D(q_{a})$ where \cite{MT2}
\begin{equation*}
D(q_{a})=\min_{\lambda \in {\Lambda}_{a}, \lambda \neq 0}|\lambda |,
\end{equation*}
is the minimal lattice distance. 
Introduce a complex sewing
parameter $\epsilon $ where $|\epsilon |\leq r_{1}r_{2}$, and excise the disk 
\begin{equation*}
\{z_{a},\left\vert z_{a}\right\vert \leq |\epsilon |/r_{\bar{a}}\}\subset 
\Sigma^{(1)}_a, 
\end{equation*}
to form a punctured torus 
\begin{equation*}
\widehat{\Sigma}^{(1)}_a=\Sigma^{(1)}_a \backslash \{z_{a},\left\vert
z_{a}\right\vert \leq |\epsilon |/r_{\bar{a}}\}.
\end{equation*}
Here and below, we use the convention 
\begin{equation*}
\overline{1}=2,\quad \overline{2}=1.  
\label{bardef}
\end{equation*}
Define the annulus
\begin{equation*}
\mathcal{A}_{a}=\{z_{a},|\epsilon |/r_{\bar{a}}\leq \left\vert
z_{a}\right\vert \leq r_{a}\}\subset \widehat{\Sigma}^{(1)}_a,
\end{equation*}%
and identify $\mathcal{A}_{1}$ and $\mathcal{A}_{2}$ as a single region $%
\mathcal{A}=\mathcal{A}_{1}\simeq \mathcal{A}_{2}$ via the sewing relation 
\begin{equation}
z_{1}z_{2}=\epsilon .  \label{pinch}
\end{equation}
In this way we obtain a compact genus two Riemann surface 
$\Sigma^{(2)}=\{ \widehat {\Sigma}^{(1)}_1
\backslash \mathcal{A}_{1} \}
\cup \{\widehat{\Sigma}^{(1)}_2 \backslash 
\mathcal{A}_{2}\}\cup \mathcal{A}$, 
parameterized by the domain \cite{MT2} 
\begin{equation}
\mathcal{D}^{\epsilon }=\{(\tau _{1},\tau _{2},\epsilon )\in \mathbb{H}_{1}
\mathbb{\times H}_{1}\mathbb{\times C}\ |\ |\epsilon |<\frac{1}{4}
D(q_{1})D(q_{2})\}.  
\label{Deps}
\end{equation}
\subsection{The Genus Two Szeg\"o Kernel in the $\epsilon$-Formalism}
\label{Subsec_eps_Szego}
On a torus the prime form is $E^{(1)}(x,y)=K^{(1)}(x-y,\tau )dx^{-\frac{1}{2}}dy^{-\frac{1}{2}}$
where 
$K^{(1)}(z,\tau )=\frac{\vartheta_{1}(z,\tau )}{\partial_{z}\vartheta _{1}(0,\tau)}$ 
and 
$\vartheta_{1}(z,\tau )
=\vartheta\left[ 
\begin{array}{c}
\frac{1}{2} 
\\
\frac{1}{2}
\end{array}
\right] (z,\tau )$ for $z\in \mathbb{C}$ and $\tau \in \mathbb{H}$. 
For $(\theta ,\phi )\neq (1,1)$ with $\theta =-\exp (-2\pi i \beta)$ and $\phi=-\exp (2\pi i \alpha)$
 the genus   one  Szeg\"{o} kernel is
\begin{equation}
S^{(1)}\left[ {
\begin{array}{c}
\theta  \\ 
\phi 
\end{array}
}\right] (x,y \vert \; \tau )=
P_{1}\left[ {\begin{array}{c}
\theta  \\ 
\phi 
\end{array}%
}\right] (x-y,\tau )dx^{\frac{1}{2}}dy^{\frac{1}{2}},
\label{S1}
\end{equation}%
where
\begin{eqnarray*}
P_{1}\left[ {%
\begin{array}{c}
\theta  \\ 
\phi 
\end{array}%
}\right] (z,\tau )&=&\frac{\vartheta^{(1)} \left[ {%
\begin{array}{c}
\alpha \\ 
\beta
\end{array}%
}\right] (z,\tau )}{\vartheta^{(1)} \left[ {%
\begin{array}{c}
\alpha 
\\ 
\beta
\end{array}
}\right] (0,\tau )}\frac{1}{K^{(1)}(z,\tau)}
\\
&=&-\sum\limits_{k\in\mathbb{Z}}\frac{q_z^{k+\lambda}}{1- \theta^{-1} q^{k+\lambda}},
\label{P1}
\end{eqnarray*}
is  a \lq twisted\rq\ Weierstrass function \cite{MTZ} and where $q_z=e^z$ and $\phi=e^{2\pi i \lambda}$ 
for $0\le \lambda<1$ (see Appendix~\ref{twisted_elliptic_functions} for details).

In \cite{TZ1} we determine the genus two Szeg\"o kernel 
\begin{equation}
S^{(2)}(x,y)=S^{(2)}\left[ {\theta}^{(2)}  \atop {\phi}^{(2)} \right] (x,y),
\label{eq:Szego2}
\end{equation}
with periodicities $\left({\theta}^{(2)}, {\phi}^{(2)} \right)= ({\theta}_{a}, {\phi}_{a})$ for $a=1,2$
on the inherited homology basis on the genus two Riemann surface $\Sigma^{(2)}$ formed by sewing
two tori $\Sigma^{(1)}_a$ in terms of genus one Szeg\"o kernel data 
$S^{(1)}_a(x,y)=S^{(1)}\left[ {\theta_a}  \atop {\phi_a} \right](x,y)$.    
Note that we exclude those Riemann theta characteristics for which $S^{(2)}$ 
exists but where one of the lower genus theta 
functions vanishes i.e. $(\theta_a,\phi_a)\neq(1,1)$ so that
$S^{(1)}_a$ exists on the torus $\Sigma^{(1)}_a$ for $a=1,2$.  
 
In \cite{TZ1} we show how to reconstruct $S^{(2)}(x,y)$ from the Laurant expansions \eqref{P1_exp} of $P_1\left[ { \theta}  \atop { \phi } \right](k,l,\tau )$ with coefficients
$C\left[ { \theta}  \atop { \phi } \right](k,l,\tau)$ and 
$D\left[ { \theta}  \atop { \phi } \right](k,l,\tau ,z)$ of \eqref{Ckldef} and \eqref{Dkldef} of 
Appendix~\ref{twisted_elliptic_functions}. In particular, we define for $k,l\ge 1$
\begin{eqnarray}
F_{a}\left[ {\theta}_{a}  \atop {\phi}_{a} \right](k,l,\tau _{a},\epsilon)
=\epsilon^{\frac{1}{2}(k+l-1)}
 C\left[ { \theta}_a  \atop { \phi }_a \right]  ( k,l,\tau _{a}). 
\label{FaCdef}
\end{eqnarray}
We let $F_a=(F_{a}\left[ {\theta}_{a}  \atop {\phi}_{a} \right](k,l,\epsilon))$ denote the infinite matrix 
indexed by $k,l\ge 1$.
We also define holomorphic $\frac{1}{2}$-forms on $\widehat{\Sigma}^{(1)}_a$ 
\begin{eqnarray}
\notag
h_{a}\left[ {\theta}_{a}  \atop {\phi}_{a} \right](k,x, \tau_a,\epsilon)
&=&
\epsilon^{\frac{k}{2} - \frac{1}{4}} \;  
D\left[ { \theta_{a} \atop \phi_{a} }\right] ( 1, k, \tau_a, x ) \; dx^{\frac{1}{2}},
\notag
\\
\bar{h}_{a}\left[ {\theta}_{a}  \atop {\phi}_{a} \right](k,y,\tau_a, \epsilon)
&=&
\epsilon^{\frac{k}{2} - \frac{1}{4}}\; 
D\left[ { \theta_{a} \atop \phi_{a} }\right] ( k, 1, \tau_a, - y) \;
 dy^{\frac{1}{2}}. 
\label{hdef}
\end{eqnarray}
We let $h_{a}(x)=(h_{a}\left[ {\theta}_{a}  \atop {\phi}_{a} \right](k,x,\tau_a, \epsilon))$ and   
$\bar{h}_{a}(y)=(\bar{h}_{a}(\left[ {\theta}_{a}  \atop {\phi}_{a} \right](k,y,\tau_a, \epsilon))$ 
denote infinite row vectors indexed by $k$. 

Recalling the $\epsilon$ sewing relation (\ref{pinch}) we note that
\begin{equation}
dz_a^{\frac{1}{2}} = (-1)^{\bar{a}}\; \xi\; \epsilon^{\frac{1}{2}}\;  \frac{dz_{\bar{a}}^{\frac{1}{2}}}
{z_{\bar{a}}},
\label{dz1dz2}
\end{equation}
where $\xi\in\{\pm \sqrt{-1}\}$ depending on the branch of the double cover of $\Sigma^{(1)}_a$ chosen.  
It is useful to introduce the infinite block matrices
\begin{equation}
\Xi=\left( 
\begin{array}{cc}
0 & \xi I 
\\ 
- \xi I & 0
\end{array}
\right),
\quad
Q=\left( 
\begin{array}{cc}
0 & \xi F_{1} 
\\ 
- \xi F_{2} & 0
\end{array}
\right), 
  \label{XFQ_def}
\end{equation}
where $I$ denotes the infinite identity matrix. Then Theorem~3.6 of \cite{TZ1} states that 
\begin{equation}
\label{s_formula}
S^{({2})}(x,y)
=\left\{ 
\begin{array}{l}
S^{(1)}_a(x,y)
 + h_{a}(x)
 \left(I-F_{\bar{a}} F_a\right)^{-1} F_{\bar{a}}  
 \bar{h}_{a}^{T}(y),  \mbox{ for }  x,y \in 
\widehat{\Sigma}^{(1)}_a, 
\\
  \xi (-1)^{\bar{a}} h_{a}(x)
 \left(I- F_{\bar{a}}F_a\right)^{-1} 
 \bar{h}_{\bar{a}}^{T}(y), \mbox{ for } 
 x\in \widehat{\Sigma}^{(1)}_a,
\; y\in \widehat{\Sigma}^{(1)}_{\bar{a}}, 
\end{array} 
\right. 
\end{equation}
where $T$ denotes the transpose. Equivalently, for 
$x,y\in \widehat{\Sigma}^{(1,1)}=\widehat{\Sigma}^{(1)}_1 \cup \widehat{\Sigma}^{(1)}_2$, 
the disconnected union of punctured tori,  
we define the forms 
\begin{eqnarray}
S^{(1,1)}(x,y)
&=&
\left\{
\begin{array}{ll}
S^{(1)}_a(x,y), & \mbox{ for } 	x,y \in \widehat{\Sigma}^{(1)}_{a}\\
0, & \mbox{ for } x \in \widehat{\Sigma}^{(1)}_{a},\ y \in \widehat{\Sigma}^{(1)}_{\bar{a}},
\end{array}
\right. \notag\\
h(x)
&=&
\left\{
\begin{array}{ll}
\left(h_{1}(x),0\right), & \mbox{ for } 	x \in \widehat{\Sigma}^{(1)}_{1}\\
\left(0,h_{2}(x)\right), & \mbox{ for } 	x \in \widehat{\Sigma}^{(1)}_{2},
\end{array}
\right.
\notag \\
\overline{h}(x)
&=&
\left\{
\begin{array}{ll}
\left(\overline{h}_{1}(x),0\right), & \mbox{ for } 	x \in \widehat{\Sigma}^{(1)}_{1}\\
\left(0,\overline{h}_{2}(x)\right), & \mbox{ for } 	x \in \widehat{\Sigma}^{(1)}_{2}.
\end{array}
\right.
\label{Shhbar}
\end{eqnarray}
Thus $h(x)$ describes an infinite row vector indexed by $k\ge 1$ and $a=1,2$ with 
$\left(h(x)\right)(k,a)=\delta_{ab}h_{b}\left[ {\theta}_{b}  \atop {\phi}_{b} \right](k,x,\tau_b, \epsilon)$ 
for $x \in \widehat{\Sigma}^{(1)}_{b}$ and similarly for $\overline{h}(x)$. 
With these definitions  \eqref{s_formula} is equivalent to
\begin{equation}
\label{sq_form}
S^{(2)}(x,y)
=S^{(1,1)}(x,y)
 + 
h(x) \Xi(I-Q)^{-1} \overline{h}^{T}(y),
\end{equation}
for $x,y \in\widehat{\Sigma}^{(1,1)}$.

Lastly, defining the determinant of $I-Q$ 
by the formal power series in $\epsilon$
\begin{equation*}
\log \det\left(I - Q\right) = Tr \log\left(I - Q\right)=-\sum_{n\ge 1}\frac{1}{n} Tr(Q^n),
\end{equation*}
it is shown in ref. \cite{TZ1} that
\begin{equation}
\label{QFdet}
\det(I-Q)=\det(I-F_1F_2), 
\end{equation}
is non-vanishing and holomorphic on $\mathcal{D}^{\epsilon }$.
\section{Vertex Operator Superalgebras} 
\label{vosas}
\subsection{General Definitions}
\label{vos}
We discuss some aspects of Vertex Operator Superalgebra theory to
establish context and notation. For more details see 
\cite{B}, \cite{FHL}, \cite{FLM}, \cite{Ka}, \cite{MN}, \cite{MT5}.
A Vertex Operator Superalgebra (VOSA) is a quadruple $(V,Y,\mathbf{1},\omega )$ 
as follows: $V$ is a superspace i.e. a complex
vector space $V=V_{\bar{0}}\oplus V_{\bar{1}}=\oplus _{\alpha }V_{\alpha }$
with index label $\alpha $ in $\mathbb{Z}/2\mathbb{Z}$ so that each $a\in V$
has a parity (fermion number) $p(a)\in \mathbb{Z}/2\mathbb{Z}$. $V$ has non-negative  
$\frac{1}{2}\mathbb{Z}$-grading with 
$$
V=\bigoplus _{r \in \frac{1}{2}\mathbb{Z}}V_{r},\mbox{ for }  \dim V_{r}<\infty,
$$  
related to the superspace grading by
\begin{equation}
V_{\bar{0}}=\bigoplus _{r \in \mathbb{Z}}V_{r},\quad V_{\bar{1}}=\bigoplus _{r \in \mathbb{Z}+\frac{1}{2}}V_{r}.\label{Vparity}
\end{equation} 
$\mathbf{1}\in V_{0}$ is the vacuum vector and $\omega
\in V_{2}$ is the conformal vector with properties described below.
$Y$ is a linear map $Y:V\rightarrow (\mathrm{End}\,V)[[z,z^{-1}]]$, for formal
variable $z$, so that for any vector $a\in V$ 
\begin{equation*}
Y(a,z)=\sum_{n\in \mathbb{Z}}a(n)z^{-n-1}.  \label{Ydefn}
\end{equation*}
The component operators (modes) $a(n)\in \mathrm{End}\,V$ are such that 
\begin{equation*}
a(n)\mathbf{1}=\delta _{n,-1}a,
\end{equation*} 
for $n\geq -1$. Furthermore, for $a\in V_{\alpha}$  
\begin{equation}
a(n):V_{\beta }\rightarrow V_{\beta+\alpha }.
\label{parityop}
\end{equation}
The vertex operators satisfy locality:
\begin{equation*}
(x-y)^{N}[Y(a,x),Y(b,y)]=0,  \label{locality}
\end{equation*}%
for all $a,b\in V$ and  $N\gg 0$, where the commutator is defined in the graded sense:  
\begin{equation}
\lbrack Y(a,x),Y(b,y)]=Y(a,x)Y(b,y)-(-1)^{p(a)p(b)}Y(b,y)Y(a,x).
\label{Ycom}
\end{equation}
The vertex operator for the vacuum is $Y(\mathbf{1},z)=Id_{V}$, whereas that
for $\omega $ is 
\begin{equation*}
Y(\omega ,z)=\sum_{n\in \mathbb{Z}}L(n)z^{-n-2},  \label{Yomega}
\end{equation*}%
where $L(n)=\omega(n+1)$ forms a Virasoro algebra for central charge $c$%
\begin{equation*}
\lbrack L(m),L(n)]=(m-n)L(m+n)+\frac{c}{12}(m^{3}-m)\delta _{m,-n}.
\label{Virasoro}
\end{equation*}%
$L(-1)$ generates translations with
\begin{equation*}
Y(L(-1)a,z)=\frac{d}{dz}Y(a,z).  \label{YL(-1)}
\end{equation*}%
$L(0)$ determines the grading with $L(0)a=wt(a)a$ for $a\in V_{r}$ and $r=wt(a)$, the weight of $a$. 
\medskip

\subsection{The Li--Zamolodchikov (Li--Z) Metric}
\label{liza}
The subalgebra $\{L(-1),L(0),L(1)\}\cong SL(2,\mathbb{C})$ associated with M\"{o}bius transformations on 
$z$ naturally acts on a VOSA (e.g. \cite{B},  \cite{Ka}). 
In particular, 
\begin{equation}
\gamma_{\lambda}=\left(
\begin{array}{cc}
0 & \lambda\\
-\lambda & 0\\	
\end{array}
\right)
:z\mapsto w=-\frac{\lambda^{2}}{z},
 \label{eq: gam_lam}
\end{equation}
is generated by 
$T_{\lambda }= \exp(\lambda L(-1))\exp(\frac{1}{\lambda}L(1))\exp(\lambda L(-1))$ where
\begin{equation}
T_{\lambda }Y(u,z)T_{\lambda }^{-1}=
Y\left(\exp(-\frac{z}{\lambda^{2}}L(1))
\left(-\frac{z}{\lambda}\right)^{-2L(0)}u,-\frac{\lambda^{2}}{z}\right).  \label{eq: Y_U}
\end{equation}
Later we will be particularly interested in the M\"{o}bius map $z\mapsto w=\epsilon/z$ associated with the 
sewing condition \eqref{pinch} with 
\begin{equation}
\lambda=-\xi\epsilon^{\frac{1}{2}},
\label{eq:lamb_eps}
\end{equation}  
with $\xi\in\{\pm \sqrt{-1}\}$ as previously introduced in \eqref{dz1dz2}.

For $u\in V$ of half-integral weight the action of $-\gamma_{\lambda}=\gamma_{-\lambda}$ is distinguished from 
that of $\gamma_{\lambda}$ whereas for integral weight they are equivalent. In particular we must distinguish the choices 
  $\lambda =\pm\sqrt{-1}$ in (\ref{eq: gam_lam}) corresponding to the
inversion map $z\mapsto z^{-1}$ normally used to define the adjoint vertex operator. Following ref. \cite{Sche} 
we therefore define 
\begin{equation}
Y^{\dagger }(u,z)=\sum_{n}u^{\dagger }(n)z^{-n-1}= T_{\lambda}Y(u,z)T_{\lambda}^{-1}. \label{eq: adj op}
\end{equation}%
One can verify that $\left(Y^{\dagger }\right)^{\dagger}(u,z)=(-1)^{2wt(u)}Y(u,z)$ for $u$ of weight $wt(u)$. 

For a quasi-primary vector $u$ (i.e. $L(1)u=0$) of weight $wt(u)$ 
\begin{equation}
u^{\dagger }(n)=\lambda^{-2wt(u)}(-\lambda^2)^{n+1}u(2wt(u)-n-2),  \label{eq: adj op qp}
\end{equation}%
e.g. $L^{\dagger}(n)=(-\lambda^2)^{n}L(-n)$. Furthermore  
\begin{equation}
Y^{\dagger }(u,w)dw^{wt(u)}= Y(u,z)dz^{wt(u)},
\label{YYform}
\end{equation} 
where for half-integral $wt(u)$ we choose the branch covering for which 
\begin{equation}
\left( \frac{dw}{dz}\right)^{wt(u)}=\left( \frac{\lambda}{z}\right)^{2wt(u)}.
\label{eq:branch}
\end{equation} 

We say a bilinear form $\langle \ ,\rangle_{\lambda}$ on $V$ is 
invariant if for all $a,b,u\in V$ \cite{Sche} 
\begin{equation}
\langle Y(u,z)a,b\rangle_{\lambda} =(-1)^{p(u)p(a)}\langle a,Y^{\dagger }(u,z)b\rangle_{\lambda}, 
\label{eq: inv bil form}
\end{equation}%
i.e. $\langle u(n)a,b\rangle_{\lambda} =(-1)^{p(u)p(a)}\langle a,u^{\dagger }(n)b\rangle_{\lambda}$. 
Thus it follows that $\langle L(0)a,b\rangle_{\lambda} =\langle a,L(0)b\rangle_{\lambda}$ so that 
$\langle a,b\rangle_{\lambda} =0$  if $wt(a)\not=wt(b)$ for homogeneous $a,b$. One also finds 
$\langle a,b\rangle_{\lambda} = \langle b,a \rangle_{\lambda}$ \cite{FHL}, \cite{Sche}. 

$\langle \ ,\rangle_{\lambda}$ is unique up to normalization if $L(1)V_{1}=V_{0}$ 
(we choose the normalization $\langle \mathbf{1} ,\mathbf{1}\rangle_{\lambda}=1$ throughout) 
and is non-degenerate if and only if $V$ is simple \cite{L}. We call such a unique non-degenerate 
symmetric bilinear form the Li--Zamolodchikov (Li--Z) metric. Given any $V$ basis $\{ u^{\alpha}\}$ we define the Li--Z dual $V$ basis $\{ \overline{u}^{\beta}\}$ where $\langle u^{\alpha} ,\overline{u}^{\beta}\rangle_{\lambda}=\delta^{\alpha\beta}$. 

\subsection{Free Fermion VOSA}
\label{Section_Rank_Two}
Consider the rank one free fermion
VOSA $V(H,\mathbb{Z}+\frac{1}{2})$ with $H=\mathbb{C}%
\psi $ for a (fermion) vector $\psi $ of parity $1$ \cite{FFR}, 
\cite{Ka} with modes obeying
\begin{equation}
\lbrack \psi (m),\psi (n)]=\psi (m)\psi (n)+\psi (n)\psi (m)=\delta
_{m+n+1,0}.  \label{Ferm_Com}
\end{equation}
The superspace is spanned by Fock vectors we denote by\footnote{Denoted by $\Psi(-\mathbf{k})$ in ref. \cite{MTZ}}
\begin{equation}
\Psi(\mathbf{k})\equiv \psi (-k_{1})\psi (-k_{2}) \ldots \psi (-k_{s})\mathbf{1},
\label{Fermion_Fock}
\end{equation}
for distinct ordered integers $1\le k_{1}<\ldots< k_{s}$ and where $\psi (k)\mathbf{1}=0$
for $k\geq 0$. The VOSA is generated by $Y(\psi, z)$ with conformal vector 
$\omega =\frac{1}{2}\psi (-2)\psi (-1)\mathbf{1}$ of central
charge $c=\frac{1}{2}$ for which $\Psi(\mathbf{k})$ has $L(0)$ 
weight $wt(\Psi(\mathbf{k}))=\sum_{1\leq i\leq s}(k_{i}-\frac{1}{2})\in \frac{1}{2}\mathbb{Z%
}$. In particular $wt(\psi )=\frac{1}{2}$. 

Since $\psi^{\dagger}(n)=\lambda^{-1} (-\lambda^2)^{n+1}\psi(-n-1)$
it follows from \eqref{eq: adj op qp} that the Fock vectors form an orthogonal basis with respect to the Li--Z metric $\langle \ ,\rangle_{\lambda}$ with
\begin{equation}
\overline{\Psi}(\mathbf{k})
=\left(-1\right )^{\left[wt(\Psi)\right]}\lambda^{2wt(\Psi)}\Psi(\mathbf{k}),
\label{eq:psidual}
\end{equation}
for $\Psi(\mathbf{k})$ of weight $wt(\Psi)$  and where $\left[x \right]$ denotes the integral part of $x$.

We next consider the rank two fermion VOSA $V(H,\mathbb{Z}+\frac{1}{2})^{\otimes 2}$, the tensor product of two copies of the rank one fermion VOSA. We employ the off-diagonal basis $\psi ^{\pm }=\frac{1}{\sqrt{2}}\left( \psi _{1}\pm i\psi _{2}\right)$ for
 fermions $\psi _{1}=\psi \otimes \mathbf{1}$ and $\psi _{2}=\mathbf{1}\otimes \psi $. 
The VOSA is generated by $Y(\psi ^{\pm },z)=\sum_{n\in \mathbb{Z} }
\psi ^{\pm }(n)z^{-n-1}$ where the modes
obey the commutation relations 
 \begin{equation*}
 \lbrack \psi ^{+}(m),\psi ^{-}(n)\rbrack=
 \delta _{m,-n-1},\quad \lbrack \psi^{+}(m),
 \psi^{+}(n)\rbrack=0,\quad \lbrack \psi^{-}(m),
 \psi^{-}(n)\rbrack=0.  \label{psiplus_minus_comm}
 \end{equation*}
The VOSA vector space $V$ is a Fock space spanned by\footnote{Denoted by $\Psi(-\mathbf{k},-\mathbf{l})$ 
in ref. \cite{MTZ}}
\begin{equation}
\Psi(\mathbf{k},\mathbf{l})\equiv \psi ^{+}(-k_{1})\ldots \psi ^{+}(-k_{s}) \psi ^{-}(-l_{1}) \ldots 
\psi^{-}(-l_{t}) \mathbf{1}, 
 \label{Fockstate}
\end{equation}
for distinct positive integers $k_{1},\ldots ,k_{s}$ and distinct $l_{1},\ldots ,l_{t}$
with $\psi^{\pm }(k) \mathbf{1}=0$ for all $k\geq 0$.
We define the conformal vector to be 
\begin{equation}
\omega =\frac{1}{2}[\psi ^{+}(-2)\psi ^{-}(-1)\mathbf{+}\psi ^{-}(-2)\psi
^{+}(-1)] \mathbf{1},  \label{omega_psi}
\end{equation}%
whose modes generate a Virasoro algebra of central charge 1. Then $\psi
^{\pm }$ has $L(0)$-weight $\frac{1}{2}$ and $\Psi(\mathbf{k},\mathbf{l})$
 has $L(0)$-weight $wt(\Psi)=\sum_{1\leq i\leq s}(k_{i}-\frac{1}{2})+\sum_{1\leq j\leq
t}(l_{j}-\frac{1}{2})$.  Similarly to \eqref{eq:psidual}, the Li--Z dual of $\Psi(\mathbf{k},\mathbf{l})$ is
\begin{equation*}
\overline{\Psi}(\mathbf{k},\mathbf{l})
=\left(-1\right )^{st}\left(-1\right )^{\left[wt(\Psi)\right]}\lambda^{2wt(\Psi)}\Psi(\mathbf{l},\mathbf{k}),
\end{equation*} 
where the $\left(-1\right )^{st}$ factor arises from the ordering chosen in \eqref{Fockstate}. For the 
parameter choice \eqref{eq:lamb_eps} we find for $\Psi(\mathbf{k},\mathbf{l})$ of parity $p_{\Psi}$ that
\begin{equation}
\overline{\Psi}(\mathbf{k},\mathbf{l})
=
\left(-1\right )^{st}(-\xi)^{p_{\Psi}}\epsilon^{wt(\Psi)}\Psi(\mathbf{l},\mathbf{k}).
\label{eq:Psibar}
\end{equation}

The weight $1$ space is $V_1=\mathbb{C}a$ for
Heisenberg  vector 
\begin{equation}
a=\psi ^{+}(-1)\psi ^{-}(-1)\mathbf{1,}  \label{adef}
\end{equation}%
with modes obeying 
\begin{equation*}
\lbrack a(m),a(n)]=m\delta_{m,-n}.
\end{equation*}%
Then $\omega=\frac{1}{2}a(-1)^{2}\mathbf{1}$ is the standard conformal vector for the Heisenberg
VOA $M$. Thus $V$ can be decomposed into irreducible $M$-modules $M\otimes e^{m}$ for $a(0)$ eigenvalue 
$m\in \mathbb{Z}$ e.g. \cite{FFR}, \cite{Ka}. Furthermore,  
$a(0)$ is a generator of continuous $V$ automorphisms $e^{2\pi i \gamma a(0)}$ for real $\gamma$. 

\section{Partition Functions and Correlation Functions on a Genus Two Riemann Surface}
\label{genustwopartfu}
In this section we consider the partition and $n$-point correlation functions for a VOSA
on a Riemann surface of genus two formed by sewing two tori. In the next section we will compute these 
quantities in the case of a rank two fermion VOSA with arbitrary automorphisms generated by $a(0)$. 

\subsection{Torus $n$-Point Correlation Functions }
\label{torusnpoint} 
We first review aspects of genus one orbifold $n$-point (correlation)
functions for twisted VOSA modules. For more details see refs. \cite{Z1}, \cite{DLM}, \cite{MT4}, \cite{DZ},  \cite{MTZ}.  

Let $\sigma \in \mathrm{Aut}(V)$ denote the parity 
(fermion number) automorphism 
\begin{equation}
\sigma a=(-1)^{p(a)}a,  \label{sigma}
\end{equation}
for all $a\in V$. Let $f,g\in \mathrm{Aut}(V)$ denote two commuting automorphisms that also commute with $\sigma$.    
Consider a $\sigma g$-twisted $V$-module $M_{\sigma g}$ 
with vertex operators $Y_{\sigma g}$ \cite{DLM}, \cite{DZ}, \cite{MTZ}.  
 We assume that $M_{\sigma g}$ is stable under $\sigma $ and $f$ i.e. both  $\sigma $ and $f$ act on $M_{\sigma g}$.
Then for vectors $v_{1},\ldots ,v_{n}\in V$ we define the torus
orbifold  $n$-point function by \cite{Z1}, \cite{MTZ}  
\begin{eqnarray}
& & Z^{(1)}\left[{f \atop g }\right](v_{1},z_{1};\ldots ;v_{n},z_{n};\tau) 
\notag
\\
& & \equiv 
\mathrm{STr}_{M_{\sigma g}}\left(
f\;Y_{\sigma g}(q_{1}^{L(0)}v_{1},q_{1})\ldots
Y_{\sigma g}(q_{n}^{L(0)}v_{n},q_{n})q^{L(0)-c/24}\right), \label{npointfunction_M_g}
\end{eqnarray}
where $q=\exp (2\pi i\tau )$, $q_{i}=\exp (z_{i})$, $i=1,\ldots, n$, for
variables $z_{1},\ldots ,z_{n}$ and where $\mathrm{STr}_{M}$ denotes the 
supertrace defined by
\begin{equation*}
\mathrm{STr}_{M}(X)=Tr_{M}(\sigma X).  
\label{Supertrace}
\end{equation*}
It follows from \eqref{parityop} that the $n$-point function \eqref{npointfunction_M_g} is non-vanishing provided
\begin{equation}
p_{1}+\ldots +p_{n}=0\mod 2,
\label{paritysum}
\end{equation} 
for parity $p_{i}=p(v_{i})$.

Taking all $v_{i}=\mathbf{1}$ in (\ref{npointfunction_M_g}) yields the 
 genus one orbifold partition function which we denote by 
$Z^{(1)}\left[{f \atop g }\right](\tau)$. Taking $n=1$ in (\ref{npointfunction_M_g}) gives the genus one 
1-point function which we denote by $Z^{(1)}\left[{f \atop g }\right] (v; \tau )$ 
and is independent of $z$.  

In order to consider modular-invariance of 
$n$-point functions at genus 1, Zhu \cite{Z1} introduced a second isomorphic
\lq square-bracket\rq\ VOSA $(V,Y[,],\mathbf{1},\tilde{\omega})$ associated to a
given VOSA $(V,Y(,),\mathbf{1},\omega )$. The new vertex operators are defined by a change of coordinates
\begin{equation*}
Y[v,z]=\sum_{n\in \mathbb{Z}}v[n]z^{-n-1}=Y(q_{z}^{L(0)}v,q_{z}-1),
\label{Ysquare}
\end{equation*}%
while the new conformal vector $\tilde{\omega}=
\omega -\frac{c}{24}\mathbf{1}$. We set $Y[\tilde{\omega},z]=\sum_{n\in \mathbb{Z}}L[n]z^{-n-2}$ 
and write $wt[v]=k$ if $L[0]v=kv$, $V_{[k]}=\{v\in V|wt[v]=k\}$. 
Only primary vectors are homogeneous with respect to both $L(0)$ and $L[0]$, 
in which case $wt(v)=wt[v]$. One can show that $n$-point functions can be expressed in terms of 1-point 
functions to find \cite{MT4}
\begin{eqnarray}
&Z^{(1)}\left[{f \atop g }\right](v_{1},z_{1};\ldots ;v_{n},z_{n};\tau)&\notag\\
&=Z^{(1)}\left[{f \atop g }\right](Y[v_{1},z_{1}-z_{n}]\ldots Y[v_{n-1},z_{n-1}-z_{n}]v_n;\tau).
&\label{Znpt1ptt}
\end{eqnarray}

\subsection{Genus Two $n$-Point Correlation Functions}
\label{genus_two_n-point}
In the $\epsilon $-sewing scheme we sew  two tori $\Sigma_{a}^{(1)}$, $a=1,2$ with modular parameters 
$\tau_{a}$ via the sewing
relation (\ref{pinch}). Similarly to ref.~\cite{MT1} for VOAs, we define the genus two orbifold
 $n$-point correlation function in the $\epsilon $-sewing scheme for a VOSA $V$ with a Li--Z metric as
 follows. Let $f_a,g_a$ be $V$ automorphisms and let $M_{\sigma g_a}$ be $\sigma g_a$-twisted $V$-modules 
stable under $\sigma$ and $f_a$ for commuting $f_a,g_a$ and $\sigma$. We combine
$f_1,g_1$ orbifold correlation functions on $\Sigma_{1}^{(1)}$ with $f_2,g_2$ orbifold correlation
functions on $\Sigma_{2}^{(1)}$.  
For $x_{1},\ldots ,x_{k}\in \Sigma_{1}^{(1)}$ with $\left\vert x_{i}\right\vert
\geq |\epsilon |/r_{2}$ and $y_{k+1},\ldots ,y_{n}\in \Sigma_{2}^{(1)}$ with $
\left\vert y_{i}\right\vert \geq |\epsilon |/r_{1}$, define the genus two orbifold $n$-point function as the following formal series in $\epsilon$
\begin{eqnarray}
& & Z^{(2)} \left[{f \atop g }\right] (v_{1},x_{1};\ldots ;  v_{k},x_{k} \vert v_{k+1},y_{k+1};
\ldots ;
v_{n},y_{n};\tau _{1},\tau _{2}, \epsilon )   
\notag
\\
\notag
\\
 \notag
& &=  \sum_{u\in V }
Z^{(1)}\left[{{f_1} \atop {g_1} }\right](Y[v_{1},x_{1}]\ldots Y[v_{k},x_{k}]u;\tau_{1}) 
\notag 
\\
\label{Z2n_pt_eps}
& &
\qquad  \cdot Z^{(1)}\left[{{f_2} \atop {g_2} }\right]
(Y[v_{k+1},y_{k+1}]\ldots Y[v_{n},y_{n}]\bar{u}; \tau _{2}), 
\end{eqnarray}
where $f$ (respectively $g$) denotes the pair $f_{1},f_{2}$ (respectively $g_{1},g_{2}$). The sum 
is taken over any $V$-basis
where $\bar{u}$ is the dual of $u$ with respect to the Li--Z metric $\langle \ ,\rangle_{\lambda}^{\mathrm{sq}}$ 
of \eqref{eq: inv bil form} as defined by the square bracket Virasoro operators $\{L[n]\}$
and with $\lambda$ of \eqref{eq:lamb_eps}.  

\begin{remark} 
\eqref{Z2n_pt_eps} reduces to the definition given in ref.~\cite{MT1} as follows. For $u,v$ of equal square 
bracket weight we have
\begin{equation}
\langle u ,v\rangle_{\lambda}^{\mathrm{sq}}=\epsilon^{-wt[u]}\langle u ,v\rangle^{\mathrm{sq}},
\label{eq:LiZcompare}
\end{equation}
where $\langle u ,v\rangle^{\mathrm{sq}}$ denotes the standard Li--Z metric corresponding to the choice $\lambda=\pm\sqrt{-1}$.
Then \eqref{Z2n_pt_eps} can be rewritten as
\begin{eqnarray*}
& & Z^{(2)} \left[{f \atop g }\right] (v_{1},x_{1};\ldots ;  v_{k},x_{k} \vert v_{k+1},y_{k+1};
\ldots ;
v_{n},y_{n};\tau _{1},\tau _{2}, \epsilon )   
\\
\\
& &=  \sum_{r\in \mathbb{Z}/2}\epsilon^{r}\sum_{u \in V_{[r]}}
Z^{(1)}\left[{{f_1} \atop {g_1} }\right](Y[v_{1},x_{1}]\ldots Y[v_{k},x_{k}]u;\tau_{1}) 
\\
& &
\qquad  \cdot Z^{(1)}\left[{{f_2} \atop {g_2} }\right]
(Y[v_{k+1},y_{k+1}]\ldots Y[v_{n},y_{n}]\bar{u}; \tau _{2}), 
\end{eqnarray*}
where here $u$ ranges over any $V_{[r]}$-basis and $\bar{u}$ is the dual of $u$ with respect to 
the standard Li--Z metric $\langle u ,v\rangle^{\mathrm{sq}}$.
\end{remark}
In the case where no states $v_{i}$ are inserted then \eqref{Z2n_pt_eps} defines the genus two partition 
(or $0$-point) function   
\begin{eqnarray}
\label{Z2_part}
Z^{(2)} \left[{f \atop g }\right] (\tau _{1},\tau _{2}, \epsilon )   
=  \sum_{u\in V }
Z^{(1)}\left[{{f_1} \atop {g_1} }\right](u;\tau_{1}) 
Z^{(1)}\left[{{f_2} \atop {g_2} }\right](\bar{u}; \tau _{2}). 
\end{eqnarray} 

The definition \eqref{Z2n_pt_eps} depends on the choice of insertion points $x_i\in \widehat{\Sigma}^{(1)}_{1}$ 
and $y_j\in \widehat{\Sigma}_{2}^{(1)}$. However, similarly to the situation for a VOA discussed in ref.~\cite{MT1},  
we may define an associated formal differential form for quasi-primary vectors as follows: 
\begin{proposition}
\label{functionformprop}
Let $v_{i}\in V$ be quasi-primary vectors of square bracket weight $wt[v_{i}]$ for $i=1,\ldots, n$.  
Let $x_{i}\in \widehat{\Sigma}^{(1)}_{1}$ and
$y_i \in \widehat{\Sigma}_{2}^{(1)}$  be related by the sewing relation 
$$x_i y_i=\epsilon=-\lambda^2.$$ 
Then the formal differential form 
\begin{eqnarray}
& & {\mathcal F}^{(2)}
\left[{f \atop g}\right]
 (v_{1}, \ldots, v_{n};\tau _{1},\tau _{2}, \epsilon )   
\notag
\\
\notag
\\
\notag
& & \equiv (-1)^{N_k}Z^{(2)} \left[{f \atop g }\right]
(v_{1},x_{1}; \ldots ;  v_{k}, x_{k}\vert  v_{k+1},y_{k+1};
\ldots ; 
v_{n}, y_{n};  \tau _{1}, \tau _{2}, \epsilon )
\\
 \label{formal}
& &  \quad
\cdot 
\prod\limits_{i=1}^{k}dx_i^{wt [v_i]}
\prod\limits_{j=k+1}^{n}dy_j^{wt [v_j]},  
\end{eqnarray}
is independent of the choice of $k=0,\ldots ,n$ where $N_{k}$ is  
the number of odd parity vectors in the set $\{v_1,\ldots ,v_{k}\}$ and where the 
 branch covering \eqref{eq:branch} is chosen with
\begin{equation*}
\left( \frac{dy_{i}}{dx_{i}}\right)^{wt [v_i]}=\left( \frac{\lambda}{x_{i}}\right)^{2wt [v_i]}.
\end{equation*} 
\end{proposition}

\textbf{Proof.} For $k\in\{1,\ldots ,n\}$ consider
\begin{eqnarray}
& & Z^{(2)} \left[{f \atop g }\right] (v_{1},x_{1};\ldots ;  v_{k},x_{k} \vert v_{k+1},y_{k+1};
\ldots ;
v_{n},y_{n})
\prod\limits_{i=1}^{k}dx_i^{wt [v_i]}
\prod\limits_{j=k+1}^{n}dy_j^{wt [v_j]}   
\notag
\\
& &=  \sum_{u\in V }
Z^{(1)}\left[{{f_1} \atop {g_1} }\right](Y[v_{1},x_{1}]\ldots 
Y[v_{k},x_{k}]u;\tau_{1})\prod\limits_{i=1}^{k}dx_i^{wt [v_i]} 
\notag
\\
& &
\qquad  \cdot Z^{(1)}\left[{{f_2} \atop {g_2} }\right]
(Y[v_{k+1},y_{k+1}]\ldots Y[v_{n},y_{n}]\bar{u}; \tau _{2})
\prod\limits_{j=k+1}^{n}dy_j^{wt [v_j]}. 
\label{Znptform}
\end{eqnarray}
We have $Y[v_{k},x_{k}]u=
\sum_{v\in V} \langle \bar{v},Y[v_{k},x_{k}]u\rangle_{\lambda}^{\mathrm{sq}} v$ where $v$ is summed over any $V$-basis. 
Since $v_{k}$ is quasi-primary,
\eqref{YYform} implies
\begin{eqnarray*}
\langle \bar{v},Y[v_{k},x_{k}]u\rangle_{\lambda}^{\mathrm{sq}}
 &=&
\langle \bar{v},Y^{\dagger}[v_{k},y_{k}]u\rangle_{\lambda}^{\mathrm{sq}}
\left( \frac{dy_{k}}{dx_{k}}\right)^{wt [v_k]}\\
&=&(-1)^{p_{k}p(v)}
\langle Y[v_{k},y_{k}]\bar{v},u\rangle_{\lambda}^{\mathrm{sq}}
\left( \frac{dy_{k}}{dx_{k}}\right)^{wt [v_k]},
\end{eqnarray*} 
using invariance \eqref{eq: inv bil form} and where $p_{k}=p(v_{k})$. Hence \eqref{Znptform} becomes
\begin{eqnarray*}
& \sum_{v\in V }(-1)^{p_{k}p(v)}
Z^{(1)}\left[{{f_1} \atop {g_1} }\right](Y[v_{1},x_{1}]\ldots 
Y[v_{k-1},x_{k-1}]v;\tau_{1})\prod\limits_{i=1}^{k-1}dx_i^{wt [v_i]}& 
\\
& 
(-1)^{p_{k}(p_{k+1}+\ldots +p_{n})}
Z^{(1)}\left[{{f_2} \atop {g_2} }\right]
(Y[v_{k},y_{k}]Y[v_{k+1},y_{k+1}]\ldots Y[v_{n},y_{n}]\bar{v}; \tau _{2})
\prod\limits_{j=k}^{n}dy_j^{wt [v_j]},& 
\end{eqnarray*}
using $\sum_{u\in V}\langle Y[v_{k},y_{k}]\bar{v},u\rangle_{\lambda}^{\mathrm{sq}}\bar{u}=Y[v_{k},y_{k}]\bar{v}$ and locality.
Finally, \eqref{paritysum} implies non-vanishing contributions arise only if $p(v)=p_{1}+\ldots+p_{k-1}$ so that \break 
$(-1)^{p_{k}p(v)}(-1)^{p_{k}(p_{k+1}+\ldots +p_{n})}=(-1)^{p_{k}}$. But $N_{k}=p_k+N_{k-1}$ where $N_{k-1}$ is the number of odd parity vectors in the set $\{v_1,\ldots, v_{k-1}\}$. Hence $(-1)^{p_{k}}=(-1)^{N_{k-1}-N_{k}}$ and the result follows.
\hfill $\square$

\medskip
  
\section{The Free Fermion VOSA} 
\label{freefermpartfu}
\subsection{Genus One}
Consider the rank 2 free fermion VOSA $V(H,\mathbb{Z}+\frac{1}{2})^{\otimes 2}$ generated by $\psi^{\pm}$. In this case, the parity automorphism \eqref{sigma} is described by $\sigma=e^{i\pi a(0)}$ for Heisenberg vector $a$. We also define two commuting automorphisms $f,g$ by\footnote{Note some notational changes from ref. \cite{MTZ}}
\begin{equation*}
\sigma f=e^{2\pi i \beta a(0)}, \quad  
\sigma g=e^{-2\pi i \alpha a(0)}, 
\end{equation*}  
for real $\alpha,\beta$. It is also convenient to define 
$\theta=-e^{-2 \pi i \beta}$, $\phi= -e^{2 \pi i \alpha}$, 
in accordance with \eqref{eq:periodicities}. The twisted partition function is then e.g. \cite{Ka}, \cite{MTZ} 
\begin{eqnarray}
Z^{(1)}
\left[{f \atop g} \right] 
\left( \tau \right)
& =& q^{\alpha^{2}/2-1/24}\prod_{l\geq 1}\left(1-\theta ^{-1}
q^{l-\frac{1}{2}+\alpha }\right)\left(1-\theta q^{l-\frac{1}{2}-\alpha }\right).  
\label{Zparth}
\end{eqnarray}
\eqref{Zparth} vanishes for $(\alpha,\beta)=(\frac{1}{2},\frac{1}{2})$ i.e. $(\theta,\phi)=(1,1)$. 
We will assume that $(\theta,\phi)\neq(1,1)$ for the remainder of this discussion. 

In ref. \cite{MTZ} it is shown by using associativity how to compute all twisted genus one $n$-point functions from a generating function which is the $2n$-point function for $n$ $\psi^{+}$ and $n$ $\psi^{-}$ vectors: 
\begin{equation}
Z^{(1)}\left[{f \atop g }\right](\psi^{+},x_{1};\psi^{-},y_{1};\ldots ;\psi^{+},x_{n};\psi^{-},y_{n};\tau)
=\det\mathbf{P} \cdot
Z^{(1)}\left[{f \atop g }\right](\tau ),  \label{FMonepointgen11}
\end{equation}%
where $\mathbf{P}$ is the $n\times n$ matrix: 
\begin{equation}
\mathbf{P}=\left( P_{1}\left[ 
{\theta  \atop \phi }
\right] (x_{i}-y_{j},\tau )\right) ,  \label{Pmatrix}
\end{equation}%
for $1\leq i,j\leq n$ and where $P_{1}\left[ 
{\theta  \atop \phi }
\right] (z,\tau)$ is the twisted Weierstrass function defined in \eqref{P1zn}. Thus, in particular, for a homogeneous square bracket weight Fock vector 
\begin{equation}
\Psi[\mathbf{k},\mathbf{l}]\equiv \psi ^{+}[-k_{1}]\ldots \psi ^{+}[-k_{s}] \psi ^{-}[-l_{1}] \ldots \psi^{-}[-l_{t}] \mathbf{1}, 
 \label{FockSq}
\end{equation}
we find that the genus one 1-point function is given by \cite{MTZ}
\begin{equation}
Z^{(1)}\left[{f \atop g }\right](\Psi[\mathbf{k},\mathbf{l}],\tau)
=\delta_{st}(-1)^{s(s-1)/2}
Z^{(1)}\left[{f \atop g }\right](\tau )\det C\left[ 
{\theta  \atop \phi }
\right]({\bf k},{\bf l},\tau ) ,  \label{FMonepointfun}
\end{equation}%
where $C\left[ 
{\theta  \atop \phi }
\right]({\bf k},{\bf l},\tau )  $ is the $s\times s$ matrix: 
\begin{equation}
C\left[ 
{\theta  \atop \phi }
\right]({\bf k},{\bf l},\tau )  =\left( C\left[ 
{\theta  \atop \phi }
\right] (k_{i},l_{j},\tau )\right),  \label{Cmatrix}
\end{equation}%
for $1\leq i,j\leq s$ as defined by \eqref{Ckldef}. Note that \eqref{FMonepointfun} is non-vanishing for $\Psi[\mathbf{k},\mathbf{l}]$ of even parity (integer weight) in agreement with  \eqref{paritysum}.

\medskip

\subsection{The Genus Two Partition Function}
We now come to the main results of this paper where for the rank two fermion VOSA we compute the genus two partition function and the generating form on the genus two Riemann surface formed by sewing together two tori as defined by \eqref{Z2n_pt_eps}. Consider commuting automorphisms $f_a,g_a$ for $a=1,2$ parameterized by 
\begin{equation*}
\sigma f_a=e^{2\pi i \beta_{a} a(0)}, \quad  
\sigma g_{a}=e^{-2\pi i \alpha_{a} a(0)}, 
\end{equation*}  
and define ${\theta_{a}}=-e^{-2 \pi i \beta_{a}}$, ${\phi}_{a}= -e^{2 \pi i \alpha_{a}}$  
where $(\theta_{a},\phi_{a})\neq (1,1)$. (The case where  $(\theta_{a},\phi_{a})=(1,1)$ will be considered elsewhere \cite{TZ4}). Recall the infinite matrices $ F_{a},Q$ of  \eqref{FaCdef} and \eqref{XFQ_def}
\begin{equation*}
F_{a}\left[ { \theta}_a  \atop { \phi }_a \right]=\left(\epsilon^{\frac{1}{2}(k+l-1)}
 C\left[ { \theta}_a  \atop { \phi }_a \right] 
 ( k,l,\tau _{a})\right),
 \quad
Q=\left( 
\begin{array}{cc}
0 & \xi F_{1}\left[ { \theta}_1  \atop { \phi }_1 \right] 
\\ 
- \xi F_{2}\left[ { \theta}_2  \atop { \phi }_2 \right] & 0
\end{array}
\right). 
\end{equation*}
We find the partition function \eqref{Z2_part} is as follows: 
\begin{theorem}
\label{Theorem_Z2_boson}
The genus two partition function for the rank two fermion VOSA is a non-vanishing holomorphic function on $\mathcal{D}^{\epsilon }$ given by 
\begin{eqnarray}
Z^{(2)}
\left[{f \atop g} \right]
\left( \tau _{1},\tau _{2}, \epsilon \right)
= 
Z^{(1)}
\left[{f_1 \atop g_1} \right]
\left( \tau_{1}\right) 
\; Z^{(1)}
\left[{f_2 \atop g_2} \right]
\left( \tau _{2}\right)\; 
\det \left( I -  Q \right). 
\label{Z2_1bos}
\end{eqnarray}
\end{theorem}

To prove this result we first note some determinant formulas for finite matrices. Let $R$ be an $N\times N$ matrix and let ${\bf k}=(k_{1},\ldots,k_{n})$ denote $n$ ordered subindices with $1\le k_{1}<\ldots <k_{n}\le N$. We refer to ${\bf k}$ as an \emph{$N$-subindex of length $n$}. Let 
\begin{equation}
R({\bf k},{\bf l})=\left( R_{k_{r}l_{s}}\right)\quad r,s=1,\ldots, n, 
\label{Rkldef}
\end{equation}
denote the $n\times n$ submatrix of $R$ indexed by a pair ${\bf k},{\bf l}$ of $N$-subindices of length $n$. We define $R({\bf k},{\bf l})=1$ in the degenerate case $n=0$.
\begin{proposition}
\label{propdet} 
Let $R$ be an $N\times N$ matrix and $I$ the identity matrix. Then 
\begin{eqnarray}
\label{detR}
\det \left(I+R\right) = 
\sum\limits_{n =0}^{N}
\sum\limits_{{\bf j}}
 \det R({\bf j},{\bf j}), 
\end{eqnarray}
where the inner sum runs over all $N$-subindices of length $n$. 
\end{proposition}
\noindent
\textbf{Proof.}
Consider $\det(I+xR)=
\sum_{\sigma\in S_{N}}\epsilon_{\sigma}
\prod_{i=1}^{N}(\delta_{i\sigma(i)}+xR_{i\sigma(i)})$
for parameter $x$ where $\epsilon_{\sigma}$ is the signature of $\sigma\in S_N$ the permutation group. Consider the subset of $S_{N}$ consisting of all permutations $\rho$ fixing \emph{at least} $N-n$
 indices. Each $\rho$ is a permutation on some 
 ${\bf j}=(j_{1},\ldots,j_{n})$, an $N$-subindex of length $n$, where the remaining $N-n$ indices are fixed. 
Then  $\det(I+xR)=\sum\limits_{0\le n \le N}a_{n}x^{n}$ for
\begin{equation*}
a_{n}=
\sum_{{\bf j}}\sum_{\rho}\epsilon_{\rho}\prod_{i=1}^{n}
R_{j_{i}\rho(j_{i})}
=\sum_{{\bf j}}\det R({\bf j},{\bf j}). 
\end{equation*}
\hfill $\square$

\begin{corollary}
\label{cordetAB} 
Let $A,B$ be $M\times M$ matrices and let $R=\left[ 
\begin{array}{cc}
0 &  tA
 \\ 
t^{-1}B  & 0
\end{array}
\right]$ be a $2M\times 2M$ block matrix for parameter $t\neq 0$. Then $\det(I+R)$ is $t$ independent and is given by
\begin{eqnarray}
\label{eqQinfty}
\det(I+R) = 
\sum\limits_{m =0}^{M} (-1)^m
\sum\limits_{{\bf k},{\bf l}}
 \det A({\bf k},{\bf l})   
 \det B({\bf l},{\bf k}), 
\end{eqnarray}
where the inner sum runs over all pairs ${\bf k},{\bf l}$ of $M$-subindices of length $m$.
\end{corollary}

\noindent
\textbf{Proof.} Clearly $I+R=
\left[ 
\begin{array}{cc}
tI_M &  0\\ 
0  & I_M
\end{array}
\right]
\left[ 
\begin{array}{cc}
I_M &  A\\ 
B  & I_M
\end{array}
\right]
\left[ 
\begin{array}{cc}
t^{-1}I_M &  0\\ 
0  & I_M
\end{array}
\right]$ 
for $M\times M$ identity matrix $I_M$ so that $\det(I+R)$ is independent of $t$. Next apply \eqref{detR} to the block matrix $R$. The block structure of $R$ and the $t$ independence of $\det(I+R)$ imply that the inner sum of \eqref{detR} runs over $2M$-indices of length $2m$ of the form ${\bf j}=(k_1,\ldots,k_m,M+l_1,\ldots,M+l_{m})$.  The pair ${\bf k},{\bf l}$ are $M$-subindices of length $m$ 
so that
\begin{equation*}
\det (I+R)=\sum_{m=0}^{M}\sum_{{\bf k}, {\bf l}}
\det\left[ 
\begin{array}{cc}
0 &  A({\bf k}, {\bf l})
 \\ 
B({\bf l},{\bf k} )  & 0
\end{array}
\right].
\end{equation*}
The result then follows.
\hfill $\square$
\medskip

\noindent
\textbf{Proof of Theorem \ref{Theorem_Z2_boson}.}  
We wish to compute the genus two partition function of \eqref{Z2_part} for the rank two fermion VOSA: 
\begin{equation*}
Z^{(2)} \left[{f \atop g }\right] (\tau _{1},\tau _{2}, \epsilon )   
=  \sum_{u\in V }
Z^{(1)}\left[{{f_1} \atop {g_1} }\right](u,\tau_{1}) 
Z^{(1)}\left[{{f_2} \atop {g_2} }\right](\bar{u}, \tau _{2}), 
\end{equation*}
where $u$ is summed over any $V$-basis and $\bar{u}$ is the square bracket Li-Z dual.
We choose the Fock basis  $\{ \Psi[{\bf k},{\bf l}] \}$ with  $1\le k_{1}< \ldots <k_{s}$ and $1\le l_{1}<\ldots < l_{m}$ of \eqref{FockSq} with square-bracket dual from \eqref{eq:Psibar}
\begin{equation}
\label{eq:Psisqbar}
\overline{\Psi}[\mathbf{k},\mathbf{l}]
=
\left(-1\right )^{sm}(-\xi)^{p_{\Psi}}\epsilon^{wt[\Psi]}\Psi[\mathbf{l},\mathbf{k}].
\end{equation} 
Furthermore, \eqref{FMonepointfun} implies the corresponding torus one point functions are non-vanishing for $m=s$ with even parity $p_{\Psi}=0$  where  
\begin{eqnarray*}
\frac{
Z^{(1)}\left[{{f_1} \atop {g_1} }\right]
(\Psi[{\bf k},{\bf l}],\tau_{1})}
{Z^{(1)}\left[{{f_1} \atop {g_1} }\right](\tau_{1})}
\notag
&=&
(-1)^{m(m-1)/2} 
\det C\left[ 
{\theta_{1}  \atop \phi_{1} }
\right]({\bf k},{\bf l},\tau_{1} ),\\
\frac{
Z^{(1)}\left[{{f_2} \atop {g_2} }\right]
(\overline{\Psi}[{\bf k},{\bf l}],\tau_{2})}
{Z^{(1)}\left[{{f_2} \atop {g_2} }\right](\tau_{2})}
&=&
(-1)^{m(m-1)/2}(-1)^m\epsilon^{wt[\Psi]}
\det C\left[ 
{\theta_{2}  \atop \phi_{2} }
\right]({\bf l},{\bf k},\tau_{2} ).
\end{eqnarray*}
Hence (suppressing the $\tau_{1},\tau_{2},\epsilon$ dependence) it follows that     
\begin{equation*}
\frac{Z^{(2)} \left[{f \atop g }\right] }
{Z^{(1)}\left[{{f_1} \atop {g_1} }\right] 
Z^{(1)}\left[{{f_2} \atop {g_2} }\right]}
=
\sum\limits_{m \ge 0}
(-1)^{m}
\sum\limits_{{\bf k},{\bf l}}\epsilon^{wt[\Psi]}
\det C\left[ 
{\theta_{1}  \atop \phi_{1} }
\right]({\bf k},{\bf l})
\det C\left[ 
{\theta_{2}  \atop \phi_{2} }
\right]({\bf l},{\bf k}). 
\end{equation*}
But $wt[\Psi]=\sum_{i=1}^{m}(k_{i}+l_{i}-1)$ so that the $\epsilon^{k_{i}+l_{j}-{\frac{1}{2}}}$ factors may be absorbed into the above $m\times m$ determinants to find 
\begin{equation*}
\frac{Z^{(2)} \left[{f \atop g }\right] }
{Z^{(1)}\left[{{f_1} \atop {g_1} }\right] 
Z^{(1)}\left[{{f_2} \atop {g_2} }\right]}
=
\sum\limits_{m \ge 0}
(-1)^{m}
\sum\limits_{{\bf k},{\bf l}}
\det F_{1}\left[ 
{\theta_{1}  \atop \phi_{1} }
\right]({\bf k},{\bf l})
\det F_{2}\left[ 
{\theta_{2}  \atop \phi_{2} }
\right]({\bf l},{\bf k}), 
\end{equation*}
with $F_{a}$ of \eqref{FaCdef}. Let $A$ and $B$ denote the finite matrices found by truncating $F_{1}$ and $F_{2}$ to an arbitrary order in $\epsilon$. 
Thus applying \eqref{eqQinfty} to $A$ and $B$ with $t=-\xi$ it follows that 
\begin{equation*}
\frac{Z^{(2)} \left[{f \atop g }\right] }
{Z^{(1)}\left[{{f_1} \atop {g_1} }\right] 
Z^{(1)}\left[{{f_2} \atop {g_2} }\right]}
=
\det(I-Q), 
\end{equation*}
as an identity between two formal series in $\epsilon$. However, it is shown in ref. \cite{TZ1} that $\det(I-Q)$ is non-vanishing and holomorphic on $\mathcal{D}^{\epsilon }$ and hence the Theorem holds.     
\hfill $\square$

\medskip

We may similarly compute the genus two partition function in the $\epsilon$-formalism for the original rank one fermion VOSA $V(H,\mathbb{Z}+\frac{1}{2})$ where, in this case, we may only construct a $\sigma$-twisted module. Then one finds:
\begin{corollary}
\label{rankonedetsq}
For the rank one free fermion VOSA $V(H,\mathbb{Z}+\frac{1}{2})$ the genus two partition function in the $\epsilon$-formalism for $f_{a},g_{a} \in \{1,\sigma  \}$ is given by 
\begin{equation}
Z^{(2)}\left[{f \atop g} \right](\tau _{1}, \tau _{2}, \epsilon)= 
Z^{(1)}\left[{f_1 \atop g_1} \right](\tau_{1}) \; 
Z^{(1)}\left[{f_2 \atop g_2} \right] (\tau_{2}) \; \det  \left(I - Q \right)^{1/2},   
\label{Z2_1bos1}
\end{equation}
where $Z^{(1)}\left[{f_{a} \atop g_{a}} \right](\tau_{a})$ is the rank one torus partition function.   
\hfill $\square$
\end{corollary}
\subsection{The Genus Two Generating Function}
\label{gen_n_point_functions}
In this section we compute the genus two generating form for all $n$-point functions for the rank two 
free fermion VOSA. This is the genus two analogue of \eqref{FMonepointgen11} and is defined by
\begin{equation}
{\cal G}_{n}^{(2)}\left[{f \atop g} \right]
(w_{1}, \ldots ,w_{n}, z_{1}, \ldots ,z_{n})
=
{\cal F}^{(2)}\left[{f \atop g} \right](\psi^{+},\psi^{-}, \ldots ,   \psi^{+},\psi^{-};  \tau_1, \tau_2, \epsilon),
\label{eq:Gentwo}
\end{equation}
the formal $2n$-form of \eqref{formal} found by alternatively inserting  $\psi^{+}$ at $w_{i}\in \widehat{\Sigma}^{(1,1)}$
 and $\psi^{-}$ at $z_{i}\in \widehat{\Sigma}^{(1,1)}$ for $i=1,\ldots, n$ where $\widehat{\Sigma}^{(1,1)}$ denotes the disconnected union of the two punctured tori. In order to describe ${\cal G}_{n}^{(2)}\left[{f \atop g} \right]$ we recall  the Szeg\"o kernels and half-forms of \eqref{eq:Szego2} and \eqref{Shhbar} and define matrices
\begin{eqnarray*}
& S^{(2)}=\left(S^{(2)}(w_{i},z_{j})\right),\quad
S^{(1,1)}=\left(S^{(1,1)}(w_{i},z_{j})\right),&\\
&H^{+}=\left((h(w_{i}))(k,a)\right), \quad
H^{-}=\left((\overline{h}(z_{i}))(l,b)\right)^{T}.&
\end{eqnarray*}
$S^{(2)}$ and  $S^{(1,1)}$ are finite matrices indexed by $w_i,z_j$ for $i,j=1,\ldots, n$; $H^{+}$  is semi-infinite with $n$ rows indexed by $w_i$ and columns
 indexed by $k\ge 1$ and $a=1,2$ and $H^{-}$  is semi-infinite with rows indexed by $l\ge 1$ and $b=1,2$ and with $n$ columns indexed by $z_j$.
We then find 
\begin{proposition}
\label{prop_Szegodet}
\begin{equation*}
\det \left[ 
\begin{array}{cc}
S^{(1,1)} & H^{+}\Xi   
\\
H^{-} & I - Q
\end{array}
\right] 
=\det  S^{(2)} \det(I - Q).
\end{equation*}
with  $Q,\Xi$ of \eqref{XFQ_def}. 
 \end{proposition}
\noindent \textbf{Proof.} 
Consider the matrix identity 
\begin{eqnarray*}
&
\left[ 
\begin{array}{cc}
S^{(1,1)} & H^{+}\Xi   
\\
H^{-} & I - Q
\end{array}
\right] 
 = &\\
&\left[ 
\begin{array}{cc}
I_{n} & H^{+}\Xi (I - Q)^{-1}  
\\
0 & I 
\end{array}
\right]
\left[ 
\begin{array}{cc}
S^{(1,1)}-H^{+}\Xi (I - Q)^{-1}H^{-}  & 0  
\\
H^{-} & I 
\end{array}
\right]
\left[ 
\begin{array}{cc}
I_{n} & 0   
\\
0 & I - Q
\end{array}
\right],&
\end{eqnarray*}
where $I_{n}$ is the $n\times n$ identity matrix. But the genus two Szeg\"o kernel of \eqref{sq_form} implies 
\begin{equation*}
\left(S^{(1,1)}-H^{+}\Xi (I - Q)^{-1}H^{-}\right)(w_{i},z_{j})= S^{(2)}(w_i, z_j). 
\end{equation*} 
The result follows on taking the determinant. \hfill $\square$
\medskip

We may next describe the generating form: 
\begin{theorem} 
\label{generating_n_point_rank_two_1}
The generating form for the rank two free fermion VOSA is given by 
\begin{eqnarray}
{\cal G}_{n}^{(2)}\left[{f \atop g} \right]
(w_{1}, \ldots ,w_{n},z_{1},\ldots ,z_{n})&=&Z^{(2)}
\left[{f \atop g} \right]
\left(\tau_{1}, \tau_2, \epsilon \right) \; \det S^{(2)}   
\label{Z2_def_epsorbi_two_point_1}. 
\end{eqnarray} 
\end{theorem}

\begin{remark}
\label{remark_Szego}
Relative to the genus two partition function, the normalized 2-point for $\psi^{+}$ and $\psi^{-}$ is given by the Szeg\"o kernel and more generally, the $2n$-point function is given by a Szeg\"o kernel determinant. This agrees with the assumed form of the higher genus fermion $2n$-point function in \cite{R} or as found by string theory methods using a Schottky parameterisation in \cite{DVPFHLS}. 
\end{remark}
In order to prove Theorem~\ref{generating_n_point_rank_two_1} we require an extension of Proposition~\ref{propdet}.
\begin{proposition}
\label{propdet2} 
Let $R$ and
$J_p=\left(
\begin{array}{cc}
0 &0\\
0 & I_{N-p}	
\end{array}
\right)$
be $N\times N$ matrices  where
$I_{N-p}$ is the identity $(N-p)\times (N-p)$ matrix for $0\le p\le N$. Then 
\begin{eqnarray}
\label{detR2}
\det \left(J_p+R\right) = 
\sum\limits_{n=0}^{N-p}
\sum\limits_{{\bf j}_p}
 \det R({\bf j}_p\, ,{\bf j}_p), 
\end{eqnarray}
where the inner sum runs over all $N$-subindices of length $n+p$ of the form 
${\bf j}_p=(1,\ldots,p,j_{1},\ldots,j_{n})$.
\end{proposition}
\noindent
\textbf{Proof.}
The proof follows along the same lines as Proposition~\ref{propdet} where here we consider $\det(J_p+xR)=x^p\sum_{\sigma\in S_{N}}\epsilon_{\sigma}
\prod_{i=1}^{p}R_{i\sigma(i)}
\prod_{i=p+1}^{N}(\delta_{i\sigma(i)}+xR_{i\sigma(i)})$.
Then $\det(J_p+xR)=x^p\sum\limits_{0\le n \le N-p}a_{n}x^{n}$ for
\begin{equation*}
a_{n}=
\sum_{{\bf j_p}}\sum_{\rho}\epsilon_{\rho}
\prod_{i=1}^{p}R_{i\rho(i)}
\prod_{r=1}^{n}R_{j_{r}\rho(j_{r})}
=\sum_{{\bf j_p}}\det R({\bf j}_p\, ,{\bf j}_p), 
\end{equation*}
where $\rho$ is a permutation of 
${\bf j}_p=(1,\ldots,p,j_{1},\ldots,j_{n})$. 
The result then follows as before.
\hfill $\square$

\begin{corollary}
\label{corTUV}
Let $A,B$ be $M\times M$ matrices and let $U$ be a $p\times M$ matrix and $W$ be a $M\times p$ matrix with $p\le M$. Define the $(p+2M)\times (p+2M)$ block matrix 
\begin{equation*}
R=\left[ 
\begin{array}{ccc}
0 & 0 & U
\\
0 & 0 & tA
 \\ 
W & t^{-1}B  & 0
\end{array}
\right], 
\end{equation*}
where $t$ is a non-zero scalar parameter.
Then $\det (J_p+R)$ is independent of $t$ and is given by
\begin{eqnarray}
\label{T}
\det (J_p+R) = 
\sum\limits_{m=p}^{M} (-1)^{m}
\sum\limits_{{\bf k},{\bf l}}
 \det U_{A}({\bf k},{\bf l})   
 \det W_{B}({\bf l},{\bf k}), 
\end{eqnarray}
where $\bf{k}$ and $\bf{l}$ are $M$-subindices of length $m-p$  and $m$ respectively. $U_{A}({\bf k},{\bf l})$ and $W_{B}({\bf l},{\bf k})$ are the $m\times m$ submatrices with components 
\begin{eqnarray*}
U_{A}({\bf k},{\bf l})_{ij}
&=&
\left \{
\begin{array}{ll}
U_{i\,l_{j}} & i=1,\ldots, p\\
A_{k_{i-p}\,l_{j}} & i=p+1,\ldots, m,	
\end{array}
\right.\\
W_{B}({\bf l},{\bf k})_{ij}
&=&
\left \{
\begin{array}{ll}
W_{l_{i}\,j} & j=1,\ldots, p\\
B_{l_{i}\,k_{j-p}} & j=p+1,\ldots, m.	
\end{array}
\right.
\end{eqnarray*}
\end{corollary}

\noindent
\textbf{Proof.} $\det (J_p+R)$ is $t$ invariant since
\begin{equation*}
(J_p+R)\vert_{t=1}=\mathrm{diag}(I_{p},t^{-1}I_{M},I_{M})\ (J_p+R)\ \mathrm{diag}(I_{p},tI_{M},I_{M}),
\end{equation*} 
for identity matrices $I_p$ and $I_M$.
$t$ invariance and the off-diagonal structure of $R$ imply that the inner sum in \eqref{detR2} is taken over $(p+2M)$-subindices of length $2m$ described by
\begin{equation*}
{\bf j}_p=(1,\ldots,p,p+k_{1},\ldots,k_{m-p},p+M+l_{1},\ldots,p+M+l_{m}),
\end{equation*} 
for $1\le k_{1}<\ldots <k_{m-p}\le M$ and $1\le l_{1}<\ldots <l_m\le M$ i.e. $\bf{k}$ and $\bf{l}$ are $M$-subindices of length $m-p$  and $m$ respectively. Hence
\begin{equation*}
\det (J_p+R)=\sum_{m=p}^{M}\sum_{{\bf k}, {\bf l}}
\det\left[ 
\begin{array}{cc}
0 &  U_{A}({\bf k}, {\bf l})
 \\ 
W_{B}({\bf l},{\bf k} )  & 0
\end{array}
\right].
\end{equation*}
The result then follows.
\hfill $\square$

\medskip

\noindent 
{\bf Proof of Theorem \ref{generating_n_point_rank_two_1}.}
Following Proposition~\ref{functionformprop} we may evaluate ${\cal G}_{n}^{(2)}\left[{f \atop g}\right]$ by inserting the quasi-primary vectors  $\psi^{\pm}$ in any way on the disconnected union of  punctured tori $\widehat{\Sigma}^{(1,1)}$. In particular, we choose $\psi^{+}$ at $w_{i}\in \widehat{\Sigma}^{(1)}_{1}$ and $\psi^{-}$ at $z_{i}\in \widehat{\Sigma}^{(1)}_{2}$ for $i=1,\ldots,n$.
 Thus, reordering operators and using \eqref{Z2n_pt_eps} and \eqref{formal} we find
\begin{eqnarray}
&&{\cal G}_{n}^{(2)}\left[{f \atop g} \right]=
{\cal G}_{n}^{(2)}\left[{f \atop g} \right]
(w_{1}, \ldots ,w_{n},z_{1},\ldots ,z_{n})\notag\\
&&=(-1)^{n(n-1)/2}(-1)^n
\sum_{u\in V }
Z^{(1)}\left[{{f_1} \atop {g_1} }\right]
(Y[\psi^{+},w_{1}]\ldots Y[\psi^{+},w_{n}]u,\tau_{1})\notag\\
&&\quad \cdot Z^{(1)}\left[{{f_2} \atop {g_2} }\right]
(Y[\psi^{-},z_{1}]\ldots Y[\psi^{-},z_{n}]\bar{u}, \tau _{2})
\prod\limits_{i=1}^{n}dw_i^{\frac{1}{2}}dz_i^{\frac{1}{2}}.
\label{GnZ1}
\end{eqnarray}
Choose the Fock basis  $\{ \Psi[{\bf k},{\bf l}] \}$ with  $1\le k_{1}< \ldots <k_{s}$ and $1\le l_{1}<\ldots < l_{m}$ of \eqref{FockSq} with square bracket dual  \eqref{eq:Psisqbar}.
The corresponding torus one point functions are non-vanishing for $n+s=m$ with parity $p_{\Psi}=n\mod 2$ from \eqref{paritysum}. Expanding \eqref{Pmatrix} using \eqref{P1_exp} one finds (see Proposition~15 of ref. \cite{MTZ} for details) 
\begin{eqnarray*}
\frac{
Z^{(1)}\left[{{f_1} \atop {g_1} }\right]
(Y[\psi^{+},w_{1}]\ldots \Psi[{\bf k},{\bf l}],\tau_{1})}
{Z^{(1)}\left[{{f_1} \atop {g_1} }\right](\tau_{1})}
&=&
(-1)^{m(m-1)/2} 
\det E_{1}({\bf k},{\bf l}),\\
\frac{Z^{(1)}\left[{{f_2} \atop {g_2} }\right]
(Y[\psi^{-},z_{1}]\ldots \overline{\Psi}[{\bf k},{\bf l}],\tau_{2})}
{Z^{(1)}\left[{{f_2} \atop {g_2} }\right](\tau_{2})}
&=&
(-1)^{m(m+1)/2} (-\xi)^{p_{\Psi}}\epsilon^{wt[\Psi]}
\det E_{2}({\bf l},{\bf k}),
\end{eqnarray*}
for $m\times m$ matrices with components
\begin{eqnarray*}
\left(E_{1}({\bf k},{\bf l})\right)_{ij}
&=&
\left \{
\begin{array}{ll}
D\left[{{\theta_{1}} \atop {\phi_{1}} }\right](1,l_{j},\tau_{1},w_{i}) & i=1,\ldots, n\\
C\left[{{\theta_{1}} \atop {\phi_{1}} }\right](k_{i},l_{j},\tau_{1}) & i=n+1,\ldots, m,	
\end{array}
\right.
\\
\left(E_{2}({\bf l},{\bf k}) \right)_{ij}
&=&
\left \{
\begin{array}{ll}
D\left[{{\theta_{2}} \atop {\phi_{2}} }\right](l_{i},1,\tau_{2},-z_{j}) & j=1,\ldots, n\\
C\left[{{\theta_{2}} \atop {\phi_{2}} }\right](l_{i},k_{j},\tau_{2}) & j=n+1,\ldots, m,
\end{array}
\right.
\end{eqnarray*}
for $C\left[{{\theta_{a}} \atop {\phi_{a}} }\right],D\left[{{\theta_{a}} \atop {\phi_{a}} }\right]$ of \eqref{Ckldef} and \eqref{Dkldef}. Since $p_{\Psi}=n\mod 2$ one finds $\xi^{p_{\Psi}}=(-1)^{n(n+1)/2}\xi^{n}$ so that altogether     
\begin{equation*}
\frac{{\cal G}_{n}^{(2)}\left[{f \atop g} \right]}
{Z^{(1)}\left[{{f_1} \atop {g_1} }\right] 
Z^{(1)}\left[{{f_2} \atop {g_2} }\right] }
=
\sum\limits_{m \ge 0}
(-1)^{m}
\sum\limits_{{\bf k},{\bf l}}\epsilon^{wt[\Psi]}\xi^{n}
\det E_{1}({\bf k},{\bf l})\det E_{2}({\bf l},{\bf k}). 
\end{equation*}
But $wt[\Psi]=\sum_{i=1}^{m-n}(k_{i}-\frac{1}{2})+\sum_{j=1}^{m}(l_{j}-\frac{1}{2})$ so that 
factors of $\epsilon^{\frac{1}{2}l_{j}-\frac{1}{4}}$ and $\epsilon^{\frac{1}{2}k_{i}-\frac{1}{4}}$ may be absorbed into the rows and columns of the above determinants. Furthermore, factors of $dw_i^{\frac{1}{2}}$ and $dz_i^{\frac{1}{2}}$ can be absorbed into the first $n$ rows and columns of $\det E_{1}$ and $\det E_{2}$ repectively. 
Lastly, a factor of $\xi$ can be absorbed into the first $n$ rows of $\det E_{1}({\bf k},{\bf l})$ to find 
\begin{equation*}
\frac{{\cal G}_{n}^{(2)}\left[{f \atop g} \right]}
{Z^{(1)}\left[{{f_1} \atop {g_1} }\right] 
Z^{(1)}\left[{{f_2} \atop {g_2} }\right] }=
\sum\limits_{m \ge 0}
(-1)^{m}
\sum\limits_{{\bf k},{\bf l}}
\det G_{1}({\bf k},{\bf l})
\det G_{2}({\bf l},{\bf k}), 
\end{equation*}
for $m\times m$ matrices
\begin{eqnarray*}
\left(G_{1}({\bf k},{\bf l})\right)_{ij}
&=&
\left \{
\begin{array}{ll}
\xi h_{1}(l_{j},\tau_{1},w_{i}) & i=1,\ldots, n\\
F_{1}(k_{i},l_{j},\tau_{1}) & i=n+1,\ldots, m,	
\end{array}
\right.
\\
\left(G_{2}({\bf l},{\bf k}) \right)_{ij}
&=&
\left \{
\begin{array}{ll}
\overline{h}_{2}
(l_{i},\tau_{2},z_{j})& j=1,\ldots, n\\
F_{2}(l_{i},k_{j},\tau_{2})
& j=n+1,\ldots, m,
\end{array}
\right.
\end{eqnarray*}
with $F_{a},h_{a}$ of \eqref{FaCdef} and \eqref{hdef}.
Finally, let $A$, $B$, $U$ and $W$ denote the finite matrices found by truncating $F_{1}$, $F_{2}$, $h_{1}(w_{i})$ and $\overline{h}_{2}
(z_{j})$ respectively to an arbitrary order in $\epsilon$. 
Thus applying Corollary~\ref{corTUV} to $A$, $B$, $U$ and $W$ with $n=p$ and $t=-\xi$ it follows that as a formal series in $\epsilon$ we have
\begin{equation*}
\frac{{\cal G}_{n}^{(2)}\left[{f \atop g} \right]}
{Z^{(1)}\left[{{f_1} \atop {g_1} }\right] 
Z^{(1)}\left[{{f_2} \atop {g_2} }\right] }
=
\det \left[ 
\begin{array}{cc}
0 & H^{+}\Xi   
\\
H^{-} & I - Q
\end{array}
\right], 
\end{equation*}
where $H^{+}\Xi=\left(0,h_{1}(l_{j},w_{i})\right)$ and $
H^{-}=\left(0,(\overline{h}_{2}(l_{i},z_{j})\right)^{T}$. Finally, using Proposition~\ref{prop_Szegodet} for $w_{i}\in \widehat{\Sigma}^{(1)}_{1}$ and $z_{i}\in \widehat{\Sigma}^{(1)}_{2}$  we find a convergent series in $\epsilon$
\begin{equation*}
\frac{{\cal G}_{n}^{(2)}\left[{f \atop g} \right]}
{Z^{(1)}\left[{{f_1} \atop {g_1} }\right] 
Z^{(1)}\left[{{f_2} \atop {g_2} }\right] }
=\det  S^{(2)} \det(I - Q),
\end{equation*}
and hence the Theorem follows on applying Theorem~\ref{Theorem_Z2_boson}. \hfill $\square$

\begin{remark} 
The other choices of the insertion points for $\psi^{\pm}$ give rise to corresponding $H^{\pm}$ and $S^{(1,1)}$ terms in Proposition~\ref{prop_Szegodet} leading to the same result \eqref{Z2_def_epsorbi_two_point_1}.  
\end{remark} 
As an illustration of the use of the generating form we compute the one-point function for the Virasoro vector $\widetilde{\omega}=
\frac{1}{2}(\psi^{+}[-2]\psi^{-}+\psi^{-}[-2]\psi^{+})$.
Let $w,z \in \widehat{\Sigma}^{(1)}_{1}$ and consider the generating form
${\cal G}_{1}^{(2)}\left[{f \atop g} \right](w,z)=S^{(2)}(w,z)  Z^{(2)}\left[{f \atop g} \right]$ (where we suppress the $\tau _{1},\tau _{2}, \epsilon$ dependence).  
Using \eqref{Znpt1ptt} we find 
\begin{eqnarray*}
&&\partial_{w}Z^{(2)} \left[{f \atop g }\right] (\psi^{+},w;\psi^{-},z)=
\partial_{w}Z^{(2)} \left[{f \atop g }\right] (Y[\psi^{+},w-z]\psi^{-},z)\\
&& =-\frac{1}{(w-z)^2}Z^{(2)} \left[{f \atop g }\right]+
Z^{(2)} \left[{f \atop g }\right] (\psi^{+}[-2]\psi^{-},z) +\ldots
\end{eqnarray*}
and similarly for $\partial_{z}$.  
Letting $S^{(2)}(w,z)=K^{(2)}(w,z)dw^{\frac{1}{2}}dz^{\frac{1}{2}}$ it follows that the Virasoro 1-point form is given by 
\begin{equation}
{\cal F}^{(2)}\left[{f \atop g }\right] (\widetilde{\omega},z)=dz^2
\lim_{w\rightarrow z}\left[
\frac{1}{2}\left(\partial_{w}-\partial_{z}\right)K^{(2)}(w,z)+\frac{1}{(w-z)^2}
\right]
Z^{(2)} \left[{f \atop g }\right].
\label{eq:Vir1pt}
\end{equation} 
An alternative expression for this is shown below in Proposition~\ref{Prop: Virasoro}.  
\medskip

\subsection{Bosonization and a Genus Two Jacobi Product Identity}
\label{bosonization} 
Consider the decomposition of the rank two fermion VOSA into irreducible modules $M\otimes e^m$ modules (for $m\in\mathbb Z$) of the Heisenberg subVOA $M$ generated by the Heisenberg state $a$.  The genus one partition function \eqref{Zparth} can thus also be expressed as (e.g. \cite{Ka}, \cite{MTZ}) 
\begin{equation*}
Z^{(1)}
\left[{f \atop g} \right] 
\left( \tau \right)
= 
\frac{e^{-2\pi i\alpha\beta}}{\eta(\tau)}
\vartheta^{(1)}\left[{\alpha \atop \beta} \right] (\tau),
\end{equation*}
for theta function \eqref{theta} and Dedekind eta-function 
$\eta (\tau )=q^{1/24}\prod\limits_{n=1}^{\infty }(1-q^{n})$. All $n$-point functions can be similarly computed in terms of Heisenberg module traces \cite{MTZ} so that the genus two partition function \eqref{Z2_1bos} can also be computed in this bosonized formalism to obtain \cite{MT1}   
\begin{equation}
Z^{(2)}
\left[{f \atop g} \right]\left(
\tau_{1},\tau_{2}, \epsilon \right)=
e^{-2\pi i\alpha\cdot\beta} Z_{M}^{(2)}(\tau_{1},\tau _{2},\epsilon )
\
\vartheta^{(2)}\left[ {\alpha \atop \beta}\right] (\Omega^{(2)} ),
\label{Z2_orbifoldpartfunc}
\end{equation}
for genus two Riemann theta function with characteristics $\alpha=(\alpha_{1},\alpha_{2}),\beta=(\beta_{1},\beta_{2})$ and where 
\begin{equation*}
Z_{M}^{(2)} (\tau_{1}, \tau_{2}, \epsilon )=
\frac{1}{\eta(\tau_1)\eta(\tau_2) \det \left(I-A_{1}A_{2}\right)^{1/2}}, 
\end{equation*}
is the genus two partition function for the rank one free Heisenberg VOA $M$. $A_{a}$ for $a=1,2$ is an infinite matrix with components indexed by $k,l\ge 1$ \cite{MT1}, \cite{MT2} 
\begin{equation*}
A_{a}(k,l,\tau_{a},\epsilon )=\epsilon^{(k+l)/2}\frac{(-1)^{k+1}(k+l-1)!}{\sqrt{kl}(k-1)!(l-1)!}%
E_{k+l}(\tau_{a}), 
\end{equation*}
for standard Eisenstein series $E_{n}(\tau)=E_{n}\left[ 
{1 \atop 1}
\right] (\tau )$. 
 Comparing with Theorem~\ref{Theorem_Z2_boson} we find a new identity relating the genus two theta function to determinants on $\mathcal{D}^{\epsilon }$ as follows
\begin{theorem}
\label{th:thetadet}
\begin{eqnarray*}
\frac{
\vartheta^{(2)}\left[ {\alpha \atop \beta}\right](\Omega^{(2)})
}
{\vartheta^{(1)} \left[ {\alpha_{1} \atop \beta_{1}}\right](\tau_1) \; 
\vartheta^{(1)} \left[ {\alpha_{2} \atop \beta_{2}}\right](\tau_2) }
=\det\left(I-A_{1}A_{2}\right)^{1/2}\; 
\det \left( I-Q\right). 
\end{eqnarray*}  
\hfill $\square$
\end{theorem}
It is shown in ref.~\cite{MT1} that $\det\left(I-A_{1}A_{2}\right)$ can be expressed as an infinite product as follows. Let $\sigma_{2n}=(k_{1},\ldots, k_{2n})$ denote a cycle permutation on $2n$ positive integers. We may canonically associate each $\sigma$ with an oriented graph $N$ consisting of $2n$ valence $2$ nodes labelled by $k_{1},\ldots,k_{2n}$. $N$ is said to be \emph{rotationless} when it admits no non-trivial rotations (a rotation being an orientation-preserving automorphism of
$N$ which preserves node labels). Lastly, we define a weight function $\zeta_{A}$ on $N$ by 
\begin{equation*}
\zeta_{A}(N)=\prod_{i=1}^{n}A_{1}(k_{2i-1},k_{2i})A_{2}(k_{2i},k_{2i+1}),
\end{equation*} 
where $k_{2n+1}\equiv k_{1}$. We then find \cite{MT1}  
\begin{equation*}
\det\left(I-A_{1}A_{2}\right)=\prod_{N\in \mathcal{R}}\left(1-\zeta_{A} (N)\right),
\end{equation*}
where $\mathcal{R}$ denotes the set of rotationless oriented cycle graphs with an even number of nodes.
This expansion can be similarly applied to $\det \left( I-Q\right)=\det \left( I-F_{1}F_{2}\right)$ with corresponding weight function $\zeta_{F} $. Hence Theorem \ref{th:thetadet} implies a genus two Jacobi product-like formula 
\begin{proposition}
\begin{eqnarray*}
\frac{
\vartheta^{(2)}\left[ {\alpha \atop \beta}\right](\Omega^{(2)})
}
{\vartheta^{(1)} \left[ {\alpha_{1} \atop \beta_{1}}\right](\tau_1) \; 
\vartheta^{(1)} \left[ {\alpha_{2} \atop \beta_{2}}\right](\tau_2) }
=\prod_{N\in \mathcal{R}}\left(1-\zeta_{A} (N)\right)^{1/2}\; 
\left(1-\zeta_{F} (N)\right).
\end{eqnarray*}  
\hfill $\square$
\end{proposition}
\begin{remark}
The bosonization procedure can be applied to obtain an alternative expression for the genus two generating form of Theorem \ref{generating_n_point_rank_two_1}
to obtain Fay's tresecant identity relating $\det S^{(2)}$ to a product of prime forms \cite{TZ4}. 
\end{remark}

\subsection{A Genus Two Ward Indentity}
We may also recompute the $1$-point function \eqref{eq:Vir1pt} for the Virasoro vector $\tilde\omega=\frac{1}{2} a[-1]a$ in the bosonized version of the rank two free fermion VOSA. 
We introduce the differential operator \cite{F1}, \cite{U}, \cite{MT1} 
\begin{equation}
\mathcal{D}=\frac{1}{2\pi i}\sum_{1\le i\le j\le 2}\nu^{(2)}_i(x)\nu^{(2)}_j(x)\frac{\partial}{%
\partial \Omega^{(2)}_{ij}},  \label{eq: Del 2form}
\end{equation}
for holomorphic 1-forms $\nu^{(2)}_i$. 
We also recall the genus two projective connection $s^{(2)}(x)$ of Appendix~\ref{PrimeForm}.
Using \eqref{Z2_orbifoldpartfunc} and results of \cite{MT1} we find 
\begin{proposition}
\label{Prop: Virasoro} 
The Virasoro 1-point form for the rank two fermion VOSA satisfies a genus two Ward identity 
\begin{equation}
\mathcal{F}^{(2)}({\tilde\omega},x;\tau_1,\tau_2,\epsilon)=
e^{-2\pi i\alpha\cdot\beta} Z_{M}^{(2)}(\tau_1,\tau_2,\epsilon) \big(\mathcal{D} + \frac{1}{12}s^{(2)}(x) %
\big )
\vartheta^{(2)}\left[ {\alpha \atop \beta}\right](\Omega^{(2)}).  \label{eq: F2omega}
\end{equation}
\end{proposition}
The Ward identity (\ref{eq: F2omega}) is similar to previous results in physics and mathematics e.g. \cite{EO}, \cite{KNTY}.

\section{Modular Invariance Properties}
\label{modular} 
We next consider the automorphic properties of the genus two partition function for the rank two fermion VOSA.  
In \cite{MTZ} we define the action of 
$\gamma=\left(
\begin{array}{cc}
a & b \\
c & d 
\end{array}
\right)
\in SL(2, \mathbb{Z})$
on a genus one orbifold partition function
$Z^{(1)}\left[{ f \atop g} \right] (\tau )$ as follows:  
\begin{equation}
\left. Z^{(1)}\left[{ f \atop g} \right] \right\vert \gamma (\tau) =Z^{(1)}\left( \gamma .\left[ { f \atop g} \right] \right) (\gamma .\tau ),  \label{ZV_gamma}
\end{equation}%
where $\gamma.\tau=\frac{a\tau +b}{c\tau+d}$ and  $\gamma .\left[ {f \atop g} \right]=\left[ {f^a g^b \atop f^cg^d} \right]$. 
For the rank two fermion VOSA we find modular invariance with \cite{MTZ}
\begin{equation}
\left. Z^{(1)}\left[{ f \atop g} \right] \right\vert \gamma (\tau) =e^{(1)}_{\gamma}\left[{ f \atop g} \right]
Z^{(1)}\left[ { f \atop g} \right]  (\tau ),
\label{eq:Z1mod}
\end{equation} 
where $e^{(1)}_{\gamma}\left[{ f \atop g} \right]\in U(1)$ is a specific multiplier system.\footnote{Note a notational change for the multiplier from that of ref. \cite{MTZ} }

In Theorem \ref{Theorem_Z2_boson} we showed that the genus two partition function is holomorphic on the domain $\mathcal{D}^{\epsilon }$ of \eqref{Deps}.  $\mathcal{D}^{\epsilon }$ is preserved under the action of $G\simeq (SL(2,\mathbb{Z})\times SL(2,\mathbb{Z}))\rtimes \mathbb{Z}_{2}$, the direct
product of the left and right torus modular groups, which are interchanged
upon conjugation by an involution $\beta$ defined as follows \cite{MT2} 
\begin{eqnarray*}
\gamma _{1}(\tau _{1},\tau _{2},\epsilon ) &=&
\left(\gamma _{1}.\tau _{1},\tau _{2},
\frac{\epsilon }{c_{1}\tau _{1}+d_{1}}\right), 
\notag 
\\\gamma _{2}(\tau _{1},\tau _{2},\epsilon ) &=&
\left(\tau _{1},\gamma _{2}.\tau _{2},\frac{\epsilon }{c_{2}\tau _{2}+d_{2}}\right), 
\notag \\
\beta (\tau _{1},\tau _{2},\epsilon ) &=&(\tau _{2},\tau _{1},\epsilon ),
\label{eq: GDeps}
\end{eqnarray*}%
for $(\gamma _{1},\gamma _{2})\in SL(2,\mathbb{Z})\times SL(2,\mathbb{Z})$
with $\gamma _{i}=\left( 
\begin{array}{cc}
a_{i} & b_{i} \\ 
c_{i} & d_{i}%
\end{array}%
\right) $.
There is a natural injection $G\rightarrow Sp (4,\mathbb{Z})$ in which the
two $SL(2,\mathbb{Z})$ subgroups are mapped to 
\begin{equation*}
\Gamma _{1}=\left\{ \left[ 
\begin{array}{cccc}
a_{1} & 0 & b_{1} & 0 \\ 
0 & 1 & 0 & 0 \\ 
c_{1} & 0 & d_{1} & 0 \\ 
0 & 0 & 0 & 1%
\end{array}%
\right] \right\} ,\quad \Gamma _{2}=\left\{ \left[ 
\begin{array}{cccc}
1 & 0 & 0 & 0 \\ 
0 & a_{2} & 0 & b_{2} \\ 
0 & 0 & 1 & 0 \\ 
0 & c_{2} & 0 & d_{2}%
\end{array}%
\right] \right\},  \label{eq: G1G2}
\end{equation*}%
and the involution is mapped to 
\begin{equation*}
\beta =\left[ 
\begin{array}{cccc}
0 & 1 & 0 & 0 \\ 
1 & 0 & 0 & 0 \\ 
0 & 0 & 0 & 1 \\ 
0 & 0 & 1 & 0%
\end{array}%
\right].  \label{eq: beta}
\end{equation*}%

In a similar way to \eqref{ZV_gamma} we define an action of 
$\gamma\in G$ on the genus two orbifold twisted partition function \eqref{Z2_part} by   
\begin{eqnarray*}
\left. Z^{(2)}\left[ {f \atop g}\right] \right \vert 
\gamma ( \tau_1,  \tau_2, \epsilon) 
= 
Z^{(2)}\left( \gamma.\left[ {f \atop g}\right] \right) \gamma.\left( \tau_1, \tau_2,\epsilon \right),  
\label{ZV2_gamma}
\end{eqnarray*}
generated by $\gamma_{i}\in \Gamma_{i}$ and $\beta$ with 
\begin{equation*}
\gamma_{1}.
\left[ 
\begin{array}{c}
f_{1}\\ f_{2} \\ g_{1} \\ g_{2}
\end{array}
\right]
=
\left[ 
\begin{array}{c}
f_{1}^{a_{1}}g_{1}^{b_{1}}\\ f_{2} \\ f_{1}^{c_{1}}g_{1}^{d_{1}} \\ g_{2}
\end{array}
\right],\quad
\gamma_{2}.
\left[ 
\begin{array}{c}
f_{1}\\ f_{2} \\ g_{1} \\ g_{2}
\end{array}
\right]
=
\left[ 
\begin{array}{c}
f_{1} \\f_{2}^{a_{2}}g_{2}^{b_{2}}\\ g_{1} \\ f_{2}^{c_{2}}g_{2}^{d_{2}}
\end{array}
\right],\quad
\beta.
\left[ 
\begin{array}{c}
f_{1}\\ f_{2} \\ g_{1} \\ g_{2}
\end{array}
\right]
=
\left[ 
\begin{array}{c}
f_{2} \\f_{1}\\ g_{2} \\ g_{1}
\end{array}
\right].
\end{equation*}

We may now describe the modular invariance of the genus two partition function for the rank two VOSA of Theorem~\ref{Theorem_Z2_boson} under the action of $G$. Define a genus two multiplier system $e^{(2)}_{\gamma}\left[{ f \atop g} \right]\in U(1)$ for $\gamma\in G$ in terms of the genus one multiplier system as follows
\begin{equation}
e^{(2)}_{\gamma_{i}}\left[{ f \atop g} \right]=
e^{(1)}_{\gamma_{i}}\left[{ f_{i} \atop g_{i}} \right],\quad
e^{(2)}_{\beta}\left[{ f \atop g} \right]=1,
\label{eq:mult2}
\end{equation}  
for $G$ generators $\gamma_{i}\in \Gamma_{i}$ and $\beta$. We then find 
\begin{theorem}
\label{two_point_modul_inv}
The genus two orbifold partition function for the rank two VOSA is modular invariant with respect to $G=(SL(2, \mathbb{Z}) \times SL(2, \mathbb{Z})) \rtimes \mathbb{Z}_2 $ with multiplier system \eqref{eq:mult2} i.e.
\begin{equation*}
\left. Z^{(2)}\left[{f \atop g}\right] \right\vert \gamma \;(\tau_1,  \tau_2, \epsilon) 
=
e^{(2)}_{\gamma} \left[{f \atop g} \right]
Z^{(2)} \left[{f \atop g}\right] 
\left(\tau_1, \tau_2, \epsilon \right).  
\end{equation*}
\end{theorem}

\noindent
{\bf Proof.}
We recall from Theorem \ref{Theorem_Z2_boson} that the genus two partition function can be expressed as
\begin{equation*}
Z^{(2)} \left[{f \atop g }\right] (\tau _{1},\tau _{2}, \epsilon )   
=  \sum_{{\bf k},{\bf l}}\left(-1\right )^{m}\epsilon^{wt[\Psi]}
Z^{(1)}\left[{{f_1} \atop {g_1} }\right]
(\Psi[{\bf k},{\bf l}],\tau_{1}) 
Z^{(1)}\left[{{f_2} \atop {g_2} }\right](\Psi[\mathbf{l},\mathbf{k}],\tau_{2}), 
\end{equation*}
for $1\le k_{1}< \ldots <k_{m}$ and $1\le l_{1}<\ldots < l_{m}$ with 
Fock basis  $\{ \Psi[{\bf k},{\bf l}] \}$ of square bracket weight $wt[\Psi]=\sum_{i=1}^{m}(k_{i}+l_{i}-1)$. Let us consider the action of $\gamma_{1}\in \Gamma_{1}$. It follows from \eqref{Ekmodular} (see also Proposition~21. of \cite{MTZ}) that
\begin{equation*}
Z^{(1)}\left(\gamma_{1}.\left[{{f_1} \atop {g_1} }\right]\right)
(\Psi[{\bf k},{\bf l}],\gamma_{1}.\tau_{1})=
e^{(1)}_{\gamma_{1}}\left[{{f_1} \atop {g_1} }\right]
(c_{1}\tau_{1}+d_{1})^{wt[\Psi]}
Z^{(1)}\left[{{f_1} \atop {g_1} }\right]
(\Psi[{\bf k},{\bf l}],\tau_{1}).
\end{equation*}
Hence from \eqref{ZV2_gamma} we find 
\begin{eqnarray*}
\left. Z^{(2)}\left[{f \atop g}\right] \right\vert \gamma_{1} 
&=&e^{(1)}_{\gamma_{1}}\left[{{f_1} \atop {g_1} }\right]
 \sum_{{\bf k},{\bf l}}\left(-1\right )^{m}
 \left(\frac{\epsilon}{c_{1}\tau_{1}+d_{1}}\right)^{wt[\Psi]}
 (c_{1}\tau_{1}+d_{1})^{wt[\Psi]}\\
 &&\cdot
Z^{(1)}\left[{{f_1} \atop {g_1} }\right]
(\Psi[{\bf k},{\bf l}],\tau_{1}) 
Z^{(1)}\left[{{f_2} \atop {g_2} }\right](\Psi[\mathbf{l},\mathbf{k}],\tau_{2})\\
&=& e^{(1)}_{\gamma_{1}}\left[{{f_1} \atop {g_1} }\right]
Z^{(2)}\left[{f \atop g}\right].  
\end{eqnarray*}
A similar result holds for $\gamma_{2}\in \Gamma_{2}$ whereas invariance under $\beta$ is obvious. The result follows. \hfill $\square$
\begin{remark}
Modular invariance can also inferred from Theorem~\ref{th:thetadet} using modular properties of the Riemann theta function together with those for the Heisenberg genus two partition function described in \cite{MT1}. 
\end{remark}
Finally, we can also obtain modular invariance for the generating form ${\cal G}_{n}^{(2)}\left[{f \atop g} \right]$ described in Theorem~\ref{generating_n_point_rank_two_1}. In particular, as is  described in \cite{TZ1}, the genus two Szeg\"o kernel of \eqref{s_formula} is invariant under the action of $G$.  Hence it follows that
\begin{theorem}
\label{n_point_modul_inv}
${\cal G}_{n}^{(2)}\left[{f \atop g} \right]$ is modular invariant with respect to $G$ with multiplier system \eqref{eq:mult2}. \hfill $\square$
\end{theorem}

\medskip

\section{Appendix}
\subsection{Some Riemann Surface Theory}
\label{PrimeForm}
Consider a compact Riemann surface $\Sigma^{(g)}$ of genus $g$ with canonical
homology cycle basis $a_{1},\ldots , a_{g}, b_{1}, \ldots , b_{g}$. In general there
exists $g$ holomorphic 1-forms $\nu^{(g)}_{i}$, $i=1,\ldots, g$ which we may
normalize  by e.g. \cite{FK}
\begin{equation*}
\oint_{a_{i}}\nu^{(g)}_{j}=2\pi i\delta_{ij}.  \label{norm}
\end{equation*}%
The genus $g$ period
matrix $\Omega^{(g)} $ is defined by 
\begin{equation*}
\Omega^{(g)}_{ij}=\frac{1}{2\pi i}\oint_{b_{i}}\nu^{(g)}_{j}, 
\label{period}
\end{equation*}%
for $i,j=1,\ldots ,g$. $\Omega^{(g)}$ is symmetric with positive imaginary part i.e. $\Omega^{(g)}\in \mathbb{H}_g$, the Siegel upper half plane.
It is useful to introduce the \emph{normalized differential of the
second kind} defined by \cite{Sp}, \cite{M}, \cite{F1}: 
\begin{equation*}
\omega^{(g)}(x,y)\sim  \frac{dxdy}{(x-y)^{2}} \quad \mbox{for } x\sim y, 
\label{omegag}
\end{equation*}
for local coordinates $x,y$, with normalization 
$
\int_{a_{i}}\omega^{(g)}(x,\cdot )=0$ for $i=1,\ldots, g$. Using the Riemann bilinear relations, one finds that 
\begin{equation*}
\nu^{(g)}_{i}(x)=\oint_{b_{i}}\omega^{(g)}(x,\cdot ). 
\end{equation*}
The projective connection $s^{(g)}$ is defined by  \cite{G} 
\begin{equation*}
s^{(g)}(x)=6\lim_{x \rightarrow y} \left (\omega^{(g)}(x,y)-\frac{dxdy}{%
(x-y)^2}\right ).  \label{eq:proj_con}
\end{equation*}
$s^{(g)}(x)$ is not a global 2-form but rather transforms under a general
conformal transformation $x\rightarrow \phi(x)$ as 
\begin{equation*}
s^{(g)}(\phi(x))=s^{(g)}(x)-\{\phi;x\}dx^2,  \label{eq: proj_con_fx}
\end{equation*}
where $\{\phi;x\}=\frac{\phi^{\prime\prime\prime}}{%
\phi^{\prime}}-\frac{3}{2}\left(\frac{\phi^{\prime\prime}}{\phi^{\prime}}%
\right)^2$ is the Schwarzian derivative. 

There exists a (nonsingular and odd) character $\left[ { \gamma} \atop {\delta} \right]$ such that 
\cite{M}, \cite{F1}  
\begin{eqnarray*}
\vartheta^{(g)} \left[ { \gamma} \atop {\delta} \right](0)=0,\qquad
 \partial _{z_i}\vartheta^{(g)} \left[ { \gamma} \atop {\delta} \right](0)\neq 0,
\label{eq:thetadelta}
\end{eqnarray*}
for the theta function with real characteristics \eqref{theta}.
Define
\begin{eqnarray*}
\zeta(x) = 
\sum_{i=1}^g 
\partial _{z_i}\vartheta^{(g)} \left[ { \gamma} \atop {\delta} \right] (0)\nu^{(g)}_i(x), 
\label{eq:zeta}
\end{eqnarray*}
a holomorphic 1-form, and let $\zeta(x)^{\frac{1}{2}}$ denote the form of weight ${\frac{1}{2}}$ on the double 
cover $\widetilde\Sigma$ of $\Sigma$. We also refer to $\zeta(x)^{\frac{1}{2}}$ as a (double-valued) ${\frac{1}{2}}$-form on $\Sigma$. 
We define the prime form $E (x, y)$ by
\begin{equation*}
\label{prime_form}
E(x,y) =\frac{ \vartheta^{(g)} \left[ {{ \gamma} \atop {\delta}} \right] 
\left( \int_{y}^{x}\nu^{(g)} \right) } 
 {\zeta(x)^{\frac{1}{2}}\zeta(y)^{\frac{1}{2}}}\sim (x-y)dx^{-\frac{1}{2}} dy^{-\frac{1}{2}} \quad \mbox{for } x\sim y,   
\end{equation*}
where $\int_{y}^{x}\nu^{(g)}= (\int_{y}^{x}\nu^{(g)}_i)\in \mathbb{C}^g$. $E(x,y)=-E(y,x)$ is 
a holomorphic differential form of weight $(-\frac{1}{2},-\frac{1}{2})$ on 
$\widetilde{\Sigma} \times \widetilde{\Sigma}$.  $E(x,y)$ has multipliers along the $a_i$ 
and $b_j$ cycles in $x$ given by $1$ and $e^{-i\pi \Omega^{(g)}_{jj}-\int_{y}^{x}\nu^{(g)}_j}$ respectively \cite{F1}. 
%
\subsection{Twisted Elliptic Functions}
\label{twisted_elliptic_functions}
Let $(\theta ,\phi )\in U(1)\times U(1)$ denote a pair of modulus one
complex parameters with $\phi =\exp (2\pi i\lambda )$ for $0\leq \lambda <1$.
 For $z\in \mathbb{C}$ and $\tau \in \mathbb{H}$ we define  \lq twisted\rq\
Weierstrass functions for $k\geq 1$ as follows \cite{MTZ}%
\begin{equation*}
P_{k}\left[ 
{\theta \atop \phi}
\right] (z,\tau )=\frac{(-1)^{k}}{(k-1)!}\sum\limits_{n\in \mathbb{Z}%
+\lambda }^{\prime }\frac{n^{k-1}q_{z}^{n}}{1-\theta ^{-1}q^{n}},
\label{Pkuv}
\end{equation*}%
for $q=q_{2\pi i\tau }$ where $\sum\limits^{\prime }$ means we omit $n=0$ if 
$(\theta ,\phi )=(1,1)$.
We  have a Laurant expansion 
\begin{equation}
P_{1}\left[{\theta \atop \phi}\right] (z,\tau )=\frac{1}{z}-\sum\limits_{n\geq 1}
E_{n}\left[{\theta \atop \phi}\right] (\tau )z^{n-1},\label{P1zn}
\end{equation}
in terms of twisted Eisenstein series for $n\geq 1$, defined
by 
\begin{eqnarray*}
E_{n}\left[ {\theta \atop \phi}\right] (\tau ) &=&-\frac{B_{n}(\lambda )}{n!}+\frac{1}{(n-1)!}%
\sum\limits_{r\geq 0}^{\prime }\frac{(r+\lambda )^{n-1}\theta
^{-1}q^{r+\lambda }}{1-\theta ^{-1}q^{r+\lambda }}  \notag \\
&& +\frac{(-1)^{n}}{(n-1)!}\sum\limits_{r\geq 1}\frac{(r-\lambda
)^{n-1}\theta q^{r-\lambda }}{1-\theta q^{r-\lambda }},  \label{Ekuv}
\end{eqnarray*}%
where $\sum\limits^{\prime }$ means we omit $r=0$ if $(\theta ,\phi )=(1,1)$
and where $B_{n}(\lambda )$ is the Bernoulli polynomial defined by%
\begin{equation*}
\frac{q_{z}^{\lambda }}{q_{z}-1}=\frac{1}{z}+\sum\limits_{n\geq 1}\frac{%
B_{n}(\lambda )}{n!}z^{n-1}.  \label{Bernoulli}
\end{equation*}

We also have Laurant expansions 
\begin{eqnarray}
& & P_{1}\left[  
{\theta \atop \phi} \right] (x-y,\tau )=\frac{1}{x-y}+\sum_{k,l\geq 1}
C\left[{\theta \atop \phi}\right] 
(k,l)\; x^{k-1}\; y^{l-1},
\notag
\\
& & P_{1}\left[ {\theta \atop \phi}\right] (z+x-y,\tau ) = 
\sum_{k,l\geq 1}D\left[{\theta \atop \phi}\right] (k,l,z)\; x^{k-1}\; y^{l-1}, 
 \label{P1_exp}
\end{eqnarray}
where for $k,l\geq 1$ we define 
\begin{eqnarray}
C\left[{\theta \atop \phi}\right] (k,l,\tau ) 
&=&(-1)^{l}\binom{k+l-2}{k-1}E_{k+l-1}\left[{\theta \atop \phi}\right] (\tau ), 
 \label{Ckldef}
\\
D\left[{\theta \atop \phi}\right] (k,l,\tau ,z) 
&=&(-1)^{k+1}\binom{k+l-2}{k-1}P_{k+l-1}\left[{\theta \atop \phi}\right] (z, \tau).
   \label{Dkldef}
\end{eqnarray}
In \cite{MTZ} we show that 
for $(\theta ,\phi )\neq (1,1)$, $E_{k}\left[{\theta \atop \phi}\right] $ is a twisted modular form of weight $k$ i.e. 
\begin{equation}
E_{k}\left( \gamma .\left[{\theta \atop \phi}\right] \right) (\gamma .\tau )
=(c\tau +d)^{k}E_{k}\left[{\theta \atop \phi}\right] (\tau ),  \label{Ekmodular}
\end{equation}
where for $\gamma=\left(\begin{array}{cc}
a &  b \\ 
c  & d 
\end{array}%
\right)\in SL(2,\mathbb{Z})$ we have $\gamma.\tau=\frac{a\tau +b}{c\tau+d}$ and  $\gamma .\left[ {\theta \atop \phi} \right]=\left[ {{\theta}^{a} {\phi}^{b} \atop {\theta}^{c}{\phi}^{d}} \right]$.

\end{document}